\documentclass{article}
\usepackage{amssymb}
\usepackage{amsfonts}
\usepackage{amsmath}

\newtheorem{theorem}{Theorem}

\newtheorem{corollary}{Corollary}

\newtheorem{lemma}{Lemma}

\newtheorem{remark}{Remark}

\numberwithin{equation}{section}
\begin{document}

\title{Reduced critical branching processes in non-favorable random
environment \thanks{%
This work was supported by the Russian Science Foundation under grant
no.24-11-00037 https://rscf.ru/en/project/24-11-00037/}}
\author{V.A.Vatutin\thanks{%
Steklov Mathematical Institute Gubkin street 8 119991 Moscow Russia Email:
vatutin@mi.ras.ru}, E.E.Dyakonova\thanks{%
Steklov Mathematical Institute Gubkin street 8 119991 Moscow Russia Email:
elena@mi-ras.ru}}
\maketitle

\begin{abstract}
Let $\left\{ Z_{n},n=0,1,2,...\right\} $ be a critical branching process in
i.i.d. random environment, $Z_{r,n}$ be the number of particles in the
process at moment $0\leq r\leq n-1$ that have a positive number of
descendants in generation $n$, and $\left\{ S_{n},n=0,1,2,...\right\} $ be
the associated random walk of $\left\{ Z_{n},n=0,1,2,...\right\} $. It is
known that if the increments of the associated random walk have zero mean
and finite variance $\sigma ^{2}$ then, for any $t\in \lbrack 0,1]$%
\begin{equation*}
\lim_{n\rightarrow \infty }\mathbf{P}\left( \frac{\log Z_{\left[ nt\right]
,n}}{\sigma \sqrt{n}}\leq x\Big|Z_{n}>0\right) =\mathbf{P}\left( \min_{t\leq
s\leq 1}B_{s}^{+}\leq x\right) ,\;x\in \lbrack 0,\infty ),
\end{equation*}%
where $\left\{ B_{t}^{+},0\leq t\leq 1\right\} $ is the Brownian meander. We
supplement this result by description of the distribution of the properly
scaled random variable $\log Z_{r,n}$ under the condition\ $\left\{
S_{n}\leq t\sqrt{k},Z_{n}>0\right\} ,$ where $t>0$ and $r,k\rightarrow
\infty $ \ in such a way that $k=o(n)$ as $n\to\infty$.

The case when the distribution of the increments of the associated random
walk belongs to the domain of attraction of a stable law is also considered.

\textbf{key words}: reduced branching processes, non-favorable random
environment, conditional limit theorem

\end{abstract}

UDK 519.218.27

\section{Introduction}

We consider a critical branching process evolving in nonfavorable random
environments. Let $\mathfrak{F}$ $=\left\{ \mathfrak{f}\right\} $ be the
space of all probability measures on $\mathbb{N}_{0}:=\{0,1,2,...\}$. For
notational reason, we identify a measure $\mathfrak{f}=\left\{ \mathfrak{f}%
(\left\{ 0\right\} ),\mathfrak{f}(\left\{ 1\right\} ),...\right\} \in $ $%
\mathfrak{F}$ with the respective probability generating function%
\begin{equation*}
f(s)=\sum_{k=0}^{\infty }\mathfrak{f}(\left\{ k\right\} )s^{k},\quad s\in
\lbrack 0,1],
\end{equation*}%
and make no difference between $\mathfrak{f}$ and $f$. Equipped with the
metric of total variation, $\mathfrak{F}$ $=\left\{ \mathfrak{f}\right\}
=\left\{ f\right\} $ becomes a Polish space. Let
\begin{equation*}
F(s)=\sum_{k=0}^{\infty }F\left( \left\{ k\right\} \right) s^{k},\quad s\in
\lbrack 0,1],
\end{equation*}%
be a random variable taking values in $\mathfrak{F}$, and let
\begin{equation*}
F_{n}(s)=\sum_{k=0}^{\infty }F_{n}\left( \left\{ k\right\} \right)
s^{k},\quad s\in \lbrack 0,1],\quad n\in \mathbb{N=}\left\{ 1,2,...\right\} ,
\end{equation*}%
be a sequence of independent probabilistic copies of $F$. The infinite
sequence $\mathcal{E}=\left\{ F_{n},n\in \mathbb{N}\right\} $ is called a
random environment.

A sequence of nonnegative random variables $\mathcal{Z}=\left\{ Z_{n},\ n\in
\mathbb{N}_{0}\right\} $ specified on\ a probability space $(\Omega ,%
\mathcal{F},\mathbf{P})$ is called a branching process in random environment
(BPRE), if $Z_{0}$ is independent of $\mathcal{E}$ and, given $\mathcal{E}$
the process $\mathcal{Z}$ is a Markov chain with
\begin{equation*}
\mathcal{L}\left( Z_{n}|Z_{n-1}=z_{n-1},\mathcal{E}=(f_{1},f_{2},...)\right)
=\mathcal{L}(\xi _{n1}+\ldots +\xi _{ny_{n-1}})
\end{equation*}%
for all $n\in \mathbb{N}$, $z_{n-1}\in \mathbb{N}_{0}$ and $%
f_{1},f_{2},...\in \mathfrak{F}$, where $\xi _{n1},\xi _{n2},\ldots $ is a
sequence of i.i.d. random variables with distribution $f_{n}.$

The sequence
\begin{equation*}
S_{0}=0,\quad S_{n}=X_{1}+...+X_{n},\ n\geq 1,
\end{equation*}%
where $X_{i}=\log F_{i}^{\prime }(1),i=1,2,...,$ is called the associated
random walk for the process $\mathcal{Z}$. The growth rate of the population
size of a BPRE essentially depends on the properties of the associated\
random walk $\mathcal{S}=\left\{ S_{n},n\geq 0\right\} $.

Let $X\overset{d}{=}X_{n},n=1,2,....$To describe the restrictions we impose
on the properties of the BPRE we introduce the set
\begin{equation*}
\mathcal{A}:=\{\alpha \in (0,2)\backslash \{1\},\,|\beta |<1\}\cup \{\alpha
=1,\beta =0\}\cup \{\alpha =2,\beta =0\}\subset \mathbb{R}^{2}.
\end{equation*}%
For a pair $(\alpha ,\beta )\in \mathcal{A}$ and a random variable $X$ we
write $X\in \mathcal{D}\left( \alpha ,\beta \right) $ if the distribution of
$X$ belongs to the domain of attraction of a stable law with density $%
g_{\alpha ,\beta }(x),x\in (-\infty ,+\infty ),$ and the characteristic
function%
\begin{equation*}
G_{\alpha ,\beta }(w)=\int_{-\infty }^{+\infty }e^{iwx}g_{\alpha .\beta
}(x)\,dx=\exp \left\{ -c|w|^{\,\alpha }\left( 1-i\beta \frac{w}{|w|}\tan
\frac{\pi \alpha }{2}\right) \right\} ,\ c>0,
\end{equation*}%
and, in addition, $\mathbf{E}X=0$ if this moment exists. Under these
assumptions there is an increasing sequence of positive numbers
\begin{equation*}
a_{n}\ =\ n^{1/\alpha }\ell (n)
\end{equation*}%
with a slowly varying at infinity sequence $\ell (1),\ell (2),\ldots ,$ such
that, as $n\rightarrow \infty $
\begin{equation*}
\left\{ \frac{S_{\left[ nt\right] }}{a_{n}},t\geq 0\right\} \Longrightarrow
\mathcal{Y}=\left\{ Y_{t},t\geq 0\right\} ,
\end{equation*}%
where $\mathcal{Y}$ is a strictly stable Levy process with the initial value
$Y_{0}=0$ and one-dimensional distributions given by the characteristic
functions%
\begin{equation*}
\mathbf{E}e^{iwY_{t}}=G_{\alpha ,\beta }(wt^{1/\alpha }),\ t\geq 0,
\end{equation*}%
and the symbol $\Longrightarrow $ stands for the weak convergence in the
spaces of c\`{a}dl\`{a}g functions endowed with Skorokhod topology. Observe
that if $X\in \mathcal{D}\left( \alpha ,\beta \right) ,$ then there exists
the limit%
\begin{equation*}
\lim_{n\rightarrow \infty }\mathbf{P}\left( S_{n}>0\right) =\rho =\mathbf{P}%
\left( Y_{1}>0\right) \in (0,1).
\end{equation*}%
Now we formulate the first condition imposed on the BPRE\ under
consideration.

\textbf{Condition B1. }The random variables\textit{\ }$X_{n},n\in \mathbb{N}%
, $\textit{\ }are independent copies of a random variable\textit{\ }$X\in
\mathcal{D}\left( \alpha ,\beta \right) $\textit{.\ }\textbf{\ }The law of%
\textit{\ }$X\in \mathcal{D}\left( \alpha ,\beta \right) $\ is absolutely
continuous and there exists\textit{\ }$n\in \mathbb{N}$\textit{\ }such%
\textit{\ }that the density $\mathbf{P}(S_{n}\in dx)/dx$\textit{\ }of\textit{%
\ }$S_{n}$\textit{\ }is bounded.\emph{\ }

\begin{remark}
\label{Rem0} The case of lattice distributions can also be analyzed by the
method described in the paper. To do this one should take care for the
correct using the properties of the renewal functions involved in the
subsequent proofs.
\end{remark}

Note that if Condition B1 is valid, then the associated random walk is
occilating and therefore BPRE under consideration is critical according to
the classification of BPRE's (see, for instance, \cite{agkv} and \cite%
{KV2017}).

Our next assumption on the environment concerns reproduction laws of
particles. Set%
\begin{equation*}
\eta =\frac{\sum_{i=1}^{\infty }i^{2}F\left( \left\{ i\right\} \right) }{%
\left( \sum_{i=0}^{\infty }iF\left( \left\{ i\right\} \right) \right) ^{2}}.
\end{equation*}

\textbf{Condition B2}. There exists $\varkappa >0$ such that\ \
\begin{equation*}
\mathbf{E}[(\log ^{+}\eta )^{\alpha +\varkappa }]\ <\ \infty ,
\end{equation*}%
where $\log ^{+}x=\log (x\vee 1)$.

This condition excludes environments with very productive particles.

Let $Z_{r,n}$ be the number of particles at moment $r\in \lbrack 0,n-1]$
that have a positive number of descendants in generation $n$ with the
additional agreement $Z_{n,n}=Z_{n}$.

\bigskip For fixed $n$ the process
\begin{equation*}
\mathcal{Z}_{red}:=\{Z_{r,n},r=0,1,...,n\}
\end{equation*}%
is called a reduced process on the interval $[0,n]$ or simply a \textit{%
reduced process}.

\bigskip It is know (see \cite{BV97} and \cite{Vat2003}) that if%
\begin{equation}
\mathbf{E}X=0,\qquad \sigma ^{2}=\mathbf{E}X^{2}\in (0,\infty ),
\label{FinVar}
\end{equation}%
and some additional technical conditions are valid, then, for any $t>0$ and $%
s\in \lbrack 0,1]$
\begin{eqnarray*}
\mathbf{P}\left( \frac{1}{\sigma \sqrt{n}}\log Z_{[sn],n}\leq
t|Z_{n}>0\right) &\sim &\mathbf{P}\left( \frac{1}{\sigma \sqrt{n}}%
\min_{[sn]\leq m\leq n}S_{m}\leq t|Z_{n}>0\right) \\
&\sim &\mathbf{P}\left( \inf_{s\leq q\leq 1}B_{q}^{+}\leq t\right)
\end{eqnarray*}%
as $n\rightarrow \infty $, where $\left\{ B_{q}^{+},0\leq q\leq 1\right\} $
is a Brownian meander, i.e., a Brownian motion conditioned to stay
nonnegative on the interval $[0,1]$.

In this paper we investigate the distribution of the random variable $%
Z_{r,n} $ in the cases when $\min (r,n-r)\rightarrow \infty $ as $%
n\rightarrow \infty ,$ and the random variable $S_{n}$ is bounded from above
by some function depending on $n$ which, in the case $\sigma ^{2}<\infty $
grows slower in order than $\sqrt{n}$. To be more specific, we study the
asymptotic behavior, as $n\rightarrow \infty $ of the conditional
distributions%
\begin{equation*}
\mathbf{P}\left( \log Z_{r,n}-S_{r}\leq ya_{m}|S_{n}\leq
ta_{k},Z_{n}>0\right)
\end{equation*}%
and
\begin{equation*}
\mathbf{P}\left( \log Z_{r,n}\leq ya_{m}|S_{n}\leq ta_{k},Z_{n}>0\right)
\end{equation*}%
depending on the ratios between the parameters $r=r(n),m=n-r$ and $k=k(n)$
as $n\rightarrow \infty .$

To formulate the main results of the paper we need some standard notation.

For two positive sequences $\left\{ c_{n},n\in \mathbb{N}\right\} $, $%
\left\{ d_{n},n\in \mathbb{N}\right\} ,$ we write as usual $c_{n}\sim d_{n}$
if $\lim_{n\rightarrow \infty }c_{n}/d_{n}=1,$ $c_{n}=O(d_{n})$ if $%
\limsup_{n\rightarrow \infty }c_{n}/d_{n}<\infty $, and $c_{n}=o(d_{n})$ or $%
c_{n}\ll d_{n}$ if $\lim_{n\rightarrow \infty }c_{n}/d_{n}=0$. Set $y\wedge
z=\min \left( y,z\right) $ and $y\vee z=\max \left( y,z\right) $.

Introduce the event%
\begin{equation*}
\mathcal{R}\left( x,n\right) =\left\{ S_{n}\leq x,Z_{n}>0\right\} .
\end{equation*}

\begin{theorem}
\label{T_reducedSmall} Assume Conditions $B1$ and $B2$. Suppose that $n\gg
k\gg m=n-r\rightarrow \infty $. Then, for any $z\in (-\infty ,+\infty )$ and
any $t>0$
\begin{equation*}
\lim_{n\rightarrow \infty }\mathbf{P}\left( \log Z_{r,n}-S_{r}\leq za_{m}|%
\mathcal{R}\left( ta_{k},n\right) \right) =\mathbf{P}\left( \min_{0\leq
s\leq 1}Y_{s}\leq z\right) .
\end{equation*}
\end{theorem}

The next theorem deals with the case $k\sim \theta m=\theta (n-r)$ for some $%
\theta \in (0,\infty )$.

\begin{theorem}
\label{T_reducedInterm} Assume that Conditions $B1$ and $B2$ are valid. If $%
k\sim \theta m=\theta m(n)\rightarrow \infty $ as $n\rightarrow \infty $ in
such a way that $m=o(n)$, then, for any $y\geq 0$ and $t>0$
\begin{equation}
\lim_{n\gg k\rightarrow \infty }\mathbf{P}\left( \log Z_{r,n}\leq ya_{m}|%
\mathcal{R}\left( ta_{k},n\right) \right) =A(\theta ^{1/\alpha }t,\theta
^{1/\alpha }\left( t\wedge y\right) ),  \label{Reduced2}
\end{equation}%
where%
\begin{equation*}
A(T,y)=\frac{\alpha \rho +1}{T^{1+\alpha \rho }}\int_{0}^{\infty }w^{\alpha
\rho }\mathbf{P}\left( -w\leq \min_{0\leq s\leq 1}Y_{s}\leq y-w;Y_{1}\leq
T-w\right) dw
\end{equation*}%
for $y\in \lbrack 0,T].$
\end{theorem}

\begin{remark}
\label{Rem1}According to Remark 4 in \cite{VD2024}
\begin{equation*}
A(T,T)=\frac{\alpha \rho +1}{T^{1+\alpha \rho }}\int_{0}^{\infty }w^{\alpha
\rho }\mathbf{P}\left( -w\leq \min_{0\leq s\leq 1}Y_{s},Y_{1}\leq T-w\right)
dw=1
\end{equation*}%
for all $T>0$. Thus, the distribution specified by the right-hand side of (%
\ref{Reduced2}) is proper.
\end{remark}

Recall that the\ meander of a strictly stable Levy process with index $%
\alpha $ is a strictly $\alpha -$stable Levy process which is assumed to be
positive on the interval $(0,1]$.

Set%
\begin{equation}
C^{\ast \ast }=\left( \int_{0}^{\infty }g^{+}\left( z\right) z^{\alpha
(1-\rho )}dz\right) ^{-1},  \label{Def_Cstar}
\end{equation}%
where $g^{+}\left( z\right) $ is the density of the distribution at moment 1
of the meander of the strictly stable\ Levy process $\mathcal{Y}$ with index
$\alpha .$

The next theorem covers some ranges in the domains $\left\{ n\gg k\gg
r\right\} ,$ $\left\{ n\gg k\sim \theta r,\theta >0\right\} $ and $\left\{
\min \left( r,n-r\right) \gg k\right\} ,$ that are not considered in
Theorems \ref{T_reducedSmall} and \ref{T_reducedInterm}, but requires rather
heavy assumptions.

\begin{theorem}
\label{T_GeneralFrac}Assume that Condition $B1$ is valid and, in addition,
there exists a constant $C\in (0,\infty )$ such that
\begin{equation*}
\mathbf{P}\left( \eta \leq C\right) =1.
\end{equation*}%
If the parameters $k$ and $r$ varies in such a way that
\begin{equation}
\log (n-r)\ll a_{k\wedge r}  \label{Cond_log}
\end{equation}%
as $n\rightarrow \infty ,$ then, for any $t>0$ and any $y\geq 0$ the
following relations are valid

1) if $n\gg k\gg r\rightarrow \infty $ then, for any $y\geq 0,$
\begin{equation*}
\lim_{r\rightarrow \infty }\mathbf{P}\left( \log Z_{r,n}\leq ya_{r}|\mathcal{%
R}\left( ta_{k},n\right) \right) =C^{\ast \ast }\mathcal{H}(y),
\end{equation*}%
where
\begin{equation*}
\mathcal{H}(y)=\int_{0}^{\infty }g^{+}\left( z\right) \left( z^{\alpha
(1-\rho )}-(z-y\wedge z)^{\alpha (1-\rho )}\right) dz;
\end{equation*}%
2) if $n\gg k\sim \theta r\rightarrow \infty $ for some $\theta >0$ then,
for any $y\geq 0,$%
\begin{equation*}
\lim_{k\rightarrow \infty }\mathbf{P}\left( \log Z_{r,n}\leq ya_{k}|\mathcal{%
R}\left( ta_{k},n\right) \right) =W(\theta ^{1/\alpha }t,\theta ^{1/\alpha
}t\wedge y),
\end{equation*}%
where
\end{theorem}

\begin{eqnarray*}
W(t,y) &=&\frac{C^{\ast \ast }(\alpha \rho +1)}{t^{\alpha \rho +1}}%
\int_{0}^{\infty }g^{+}\left( z\right) dz\int_{(z-y)\vee 0}^{z}q^{\alpha
\rho }(t-z+q)^{\alpha (1-\rho )}dq \\
&&\qquad +\frac{C^{\ast \ast }\alpha (1-\rho )}{t^{\alpha \rho +1}}%
\int_{0}^{\infty }g^{+}\left( z\right) dz\int_{(z-y)\vee 0}^{z}\left(
t-z+q\right) ^{\alpha \rho +1}q^{\alpha (1-\rho )-1}dq;
\end{eqnarray*}%
3) if $\min \left( r,n-r\right) \gg k\rightarrow \infty $ then, for any $%
y\geq 0,$%
\begin{equation*}
\lim_{k\rightarrow \infty }\mathbf{P}\left( \log Z_{r,n}\leq ya_{k}|\mathcal{%
R}\left( ta_{k},n\right) \right) =1-\left( 1-\frac{t\wedge y}{t}\right)
^{\alpha \rho +1}.
\end{equation*}

\begin{remark}
\label{Rem3} It is not difficult to check that $\eta =1$ almost surely (a.s)
if%
\begin{equation*}
F(s)=e^{\lambda (1-s)},\;\lambda >0,
\end{equation*}%
and $\eta =2$ a.s if
\begin{equation*}
F(s)=\frac{q}{1-ps},\;p+q=1,\;pq>0.
\end{equation*}
\end{remark}

\begin{remark}
\label{Rem3a} If $X_{i}\in \mathcal{D}\left( \alpha ,\beta \right)
,i=1,2,...,$ and \ either $n\gg k\gg r$ or $n\gg k\sim \theta r,$ then
condition (\ref{Cond_log}) is equivalent to the assumption $\log n\ll
a_{k\wedge r}$. Note also that the restriction $\min \left( r,n-r\right) \gg
k$ includes the cases when $r\sim cn$ for some $c\in (0,1).$ Under this
assumption and the validity of (\ref{FinVar}) relation (\ref{Cond_log})
looks as follows
\begin{equation*}
\log ^{2}n\ll k\ll n.
\end{equation*}
\end{remark}

\section{Some properties of random walk}

In the sequel we denote by $C_{1},C_{2},...$ some absolute constants. We
also agree that one and the same constant may have different values in
different formulas.

We recall that a positive sequence $\left\{ c_{n},n\in \mathbb{N}\right\} $
(or a real function $c(x),x\geq 0)$ is said to be regularly varying at
infinity with index $\gamma \in \mathbb{R}$ , denoted $\left\{ c_{n},n\in
\mathbb{N}\right\} \in R_{\gamma }$ ($c(x)\in R_{\gamma }$) if $c_{n}\sim
n^{\gamma }l(n)$ as $n\rightarrow \infty $ $(c(x)\sim x^{\gamma }l(x)$ as $%
x\rightarrow \infty )$, where $l(x)$ is a slowly varying function.

Set $S_{0}:=0,$ $\tau _{0}^{\pm }:=0,$ and denote%
\begin{equation*}
\tau _{k}^{-}:=\inf \left\{ n>\tau _{k-1}^{-}:S_{n}\leq S_{\tau
_{k-1}^{-}}\right\} ,\;k\geq 1,
\end{equation*}%
for weak descending ladder variables and%
\begin{equation*}
\tau _{k}^{+}:=\inf \left\{ n>\tau _{k-1}^{+}:S_{n}\geq S_{\tau
_{k-1}^{+}}\right\} ,\;k\geq 1,
\end{equation*}%
for weak ascending ladder variables. \ Put
\begin{equation*}
H_{k}^{\pm }:=\pm S_{\tau _{k}^{\pm }}
\end{equation*}%
and denote%
\begin{equation*}
\zeta =\mathbf{P}\left( H_{1}^{+}=0\right) =\mathbf{P}\left(
H_{1}^{-}=0\right) \in \lbrack 0,1).
\end{equation*}

For $x\geq 0$ introduce renewal functions%
\begin{equation*}
V^{\pm }(x)=\sum_{k=0}^{\infty }\mathbf{P}\left( H_{k}^{\pm }\leq x\right)
=\sum_{k=0}^{\infty }\sum_{n=0}^{\infty }\mathbf{P}\left( \tau _{k}^{\pm
}=n,\pm S_{n}\leq x\right) .\quad
\end{equation*}%
Observe that $V^{\pm }(x)$ are non-decreasing and%
\begin{equation}
V^{\pm }(0)=\sum_{k=0}^{\infty }\mathbf{P}\left( H_{k}^{\pm }=0\right) =%
\frac{1}{1-\zeta }.  \label{V_zero}
\end{equation}%
It is known (see, for instance, \cite{Rog1971}, \cite{Sin57}) that if $X\in
\mathcal{D}\left( \alpha ,\beta \right) $, then
\begin{equation}
\mathbf{P}\left( \tau _{1}^{+}>n\right) \in R_{-\rho },\quad V^{+}(x)\in
R_{\alpha \rho },  \label{Regular1}
\end{equation}%
and%
\begin{equation}
\mathbf{P}\left( \tau _{1}^{-}>n\right) \in R_{-(1-\rho )},\quad V^{-}(x)\in
R_{\alpha (1-\rho )}.  \label{Regular2}
\end{equation}%
It follows from (\ref{Regular1}) that, as $x\rightarrow \infty $
\begin{equation}
\left( \alpha \rho +1\right) \int_{0}^{x}V^{+}(w)dw\sim xV^{+}(x).
\label{AsympV}
\end{equation}

We also introduce for $x\geq 0$ renewal functions
\begin{equation*}
\underline{{V}}^{\pm }(x):=\sum_{k=0}^{\infty }\mathbf{P}\left( H_{k}^{\pm
}<x\right) =\sum_{k=0}^{\infty }\sum_{n=0}^{\infty }\mathbf{P}\left( \tau
_{k}^{\pm }=n,\pm S_{n}<x\right) .
\end{equation*}%
If the distribution of $X_{1}$ is absolutely continuous, then%
\begin{equation*}
\underline{{V}}^{\pm }(x)=V^{\pm }(x),\;x\geq 0.
\end{equation*}

The renewal functions we just have defined are constructed by weak ladder
moments and ladder heights. Strict ladder moments $\left\{ \hat{\tau}%
_{k}^{\pm },k\geq 0\right\} $ and strict ladder heights $\left\{ \hat{H}%
_{k}^{\pm },k\geq 0\right\} $ defined as $\hat{\tau}_{0}^{\pm }:=0,$\ $\hat{H%
}_{0}^{\pm }:=0$ and, for $k\geq 1$
\begin{equation*}
\hat{\tau}_{k}^{\pm }:=\inf \left\{ n>\hat{\tau}_{k-1}^{\pm }:\pm S_{n}>\pm
S_{\hat{\tau}_{k-1}^{\pm }}\right\} ,\quad \hat{H}_{k}^{\pm }:=\pm S_{\tau
_{k}^{\pm }},
\end{equation*}%
are also often used in studying properties of random walks and branching
processes in random environment (see, for instance, \cite{agkv}, \cite%
{ABGV2011},\cite{Don12}, \cite{VW09},\cite{VD2022}).

The sequences $\left\{ \hat{H}_{k}^{\pm },k\geq 0\right\} $ generate the
renewal functions
\begin{eqnarray*}
\hat{V}^{\pm }(x) &:&=\sum_{k=0}^{\infty }\mathbf{P}\left( \hat{H}_{k}^{\pm
}\leq x\right) =\sum_{k=0}^{\infty }\sum_{n=0}^{\infty }\mathbf{P}\left(
\hat{\tau}_{k}^{\pm }=n,\pm S_{n}\leq x\right) , \\
\underline{\hat{V}}^{\pm }(x) &:&=\sum_{k=0}^{\infty }\mathbf{P}\left( \hat{H%
}_{k}^{\pm }<x\right) =\sum_{k=0}^{\infty }\sum_{n=0}^{\infty }\mathbf{P}%
\left( \hat{\tau}_{k}^{\pm }=n,\pm S_{n}<x\right) .
\end{eqnarray*}%
It is known (see \cite[Section XII.1, equation (1.13) ]{Fel}) that%
\begin{equation}
\hat{V}^{\pm }(x)=\left( 1-\zeta \right) V^{\pm }(x),\quad \underline{\hat{V}%
}^{\pm }(x)=\left( 1-\zeta \right) \underline{{V}}^{\pm }(x).
\label{RenewRelation}
\end{equation}

Denote%
\begin{equation}
b_{n}=\frac{1}{na_{n}},\;n\geq 1,  \label{Defb}
\end{equation}%
and set%
\begin{eqnarray*}
L_{r,n} &:&=\min_{r\leq i\leq n}S_{i},\quad L_{n}:=L_{0,n},\quad
M_{n}:=\max_{1\leq i\leq n}S_{i}, \\
\tau _{r,n} &:&=\min \left\{ r\leq i\leq n:S_{i}=L_{r,n}\right\} ,\quad \\
\tau _{n} &:&=\tau _{0,n}=\min \left\{ 0\leq i\leq n:S_{i}=L_{n}\right\} .
\end{eqnarray*}%
Introduce the event%
\begin{equation*}
\mathcal{B}(x,n):=\{S_{n}\leq x,L_{n}\geq 0\}
\end{equation*}%
which often appears below in our arguments.

Sometimes we will consider random walk under $\mathbf{P}$ \ started from an
arbitrary point $w\in \mathbb{R}$. In this case we will use the notation $%
\mathbf{P}_{w}$ \ with the natural agreement $\mathbf{P}=\mathbf{P}_{0}$.
According to Corollary 2 in \cite{VD2022} (written in the notation and the
assumptions of the present paper)
\begin{equation}
\mathbf{P}_{w}\left( \mathcal{B}(x,n)\right) \sim g_{\alpha ,\beta
}(0)V^{-}(w)b_{n}\int_{0}^{x}V^{+}(u)du  \label{BasicAsymptotic}
\end{equation}%
as $n\rightarrow \infty $ uniformly in $w,x\geq 0$ such that $\max (x,w)\in
(0,\delta _{n}a_{n}],$ where $\delta _{n}\rightarrow 0$ as $n\rightarrow
\infty $. Moreover, there is a constant $C$ such that%
\begin{equation*}
\mathbf{P}_{w}\left( x<S_{n}\leq x+1,L_{n}\geq 0\right) \leq
CV^{-}(w)b_{n}V^{+}(x)
\end{equation*}%
and%
\begin{equation}
\mathbf{P}_{w}\left( \mathcal{B}(x,n)\right) \leq
CV^{-}(w)b_{n}\int_{0}^{x}V^{+}(u)du  \label{Rough1}
\end{equation}%
for all $w,x\geq 0$ and all $n\geq 1$.

The distribution of the random variable $S_{\tau _{r,n}}$ under the
condition $\mathcal{B}(x,n)$ plays an important role in studying the
distribution of the random variable $Z_{r,n}$. For this reason and
convenience of references we below formulate a number of theorems
established in \cite{VDD2024} and describing the limiting behavior of the
conditional probabilities%
\begin{equation*}
\mathbf{P}_{w}\left( S_{\tau _{r,n}}\leq ya_{r}|\mathcal{B}(ta_{k},n)\right)
\text{ \ \ and \ \ }\mathbf{P}_{w}\left( S_{\tau _{r,n}}\leq ya_{k}|\mathcal{%
B}(ta_{k},n)\right) \text{\ },\quad t>0,
\end{equation*}%
under various assumptions concerning relations between the parameters $r$
and $k$ as $n\rightarrow \infty $.\bigskip

\begin{theorem}
\label{T_beginningNew} If Condition B1 is valid, then for any fixed $w\geq
0, $ $y\geq 0$ and $t>0$ the following relations are valid:

1) (\cite[Theorem 1]{VD2024})
\begin{equation*}
\lim_{n\gg k\gg r\rightarrow \infty }\mathbf{P}_{w}\left( S_{\tau
_{r,n}}\leq ya_{r}|\mathcal{B}(ta_{k},n)\right) =C^{\ast \ast }\mathcal{H}%
(y),
\end{equation*}%
where $C^{\ast \ast }$ is the same as in (\ref{Def_Cstar})

2) if $k\sim \theta r,\theta \in (0,\infty )$, then (\cite[Theorem 2]{VD2024}%
)
\begin{equation*}
\lim_{n\gg k\rightarrow \infty }\mathbf{P}_{w}\left( S_{\tau _{r,n}}\leq
ya_{r}|\mathcal{B}(ta_{k},n)\right) =W(t\theta ^{1/\alpha },\theta
^{1/\alpha }t\wedge y);
\end{equation*}

3) (\cite[Theorem 3]{VD2024})
\begin{equation*}
\lim_{\min (r,n-r)\gg k\rightarrow \infty }\mathbf{P}_{w}\left( S_{\tau
_{r,n}}\leq ya_{k}|\mathcal{B}(ta_{k},n)\right) =1-\left( 1-\frac{t\wedge y}{%
t}\right) ^{\alpha \rho +1};
\end{equation*}

4) if $k\sim \theta (n-r),\theta \in (0,\infty ),$ and $n-r=o(n)$ as $%
n\rightarrow \infty $, then \cite[Theorem 4]{VD2024}
\begin{equation*}
\lim_{n\gg k\rightarrow \infty }\mathbf{P}_{w}\left( S_{\tau _{r,n}}\leq
ya_{k}|\mathcal{B}(ta_{k},n)\right) =A(\theta ^{1/\alpha }t,\theta
^{1/\alpha }\left( t\wedge y\right) ).
\end{equation*}%
\bigskip
\end{theorem}

The next theorem has a slightly different form in comparison with the
statements of Theorem \ref{T_beginningNew}.

\begin{theorem}
\label{T_Minim}(\cite[Theorem 5]{VD2024}) If Condition B1 is valid and $%
m=n-r $, then, for any fixed $w\geq 0$ and $t>0$
\begin{equation*}
\lim_{n\gg k\gg m\rightarrow \infty }\mathbf{P}_{w}\left( S_{\tau
_{r,n}}-S_{r}\leq ya_{m}|\mathcal{B}(ta_{k},n)\right) =\mathbf{P}\left(
\min_{0\leq s\leq 1}Y_{s}\leq y\right) .
\end{equation*}
\end{theorem}

\section{\protect\bigskip Conditional limit theorem for random walks}

In this section we will forget about branching processes in random
environment for a while \ and assume, keeping the previous notation, that we
have a random walk $\mathcal{S}$, generated by a sequence of i.i.d. random
variables. To analyse properties of such random walk we need consider\
conditional laws $\mathbf{P}_{x}^{+}\left( \cdot \right) $ and $\mathbf{P}%
_{-x}^{-}\left( \cdot \right) ,x\geq 0$, the laws of the random walk under $%
\mathbf{P}$ started at $\pm x$ at time $n=0$ and conditioned to stay
non-negative or negative for all times (see, for example, \cite{BD94}, \cite%
{CC2008} and \cite{VW09}). These laws are specified for all $x\geq 0$, $n\in
\mathbb{N}$, and sets $B_{n}$ belonging to the $\sigma$-algebra generated
by the random variables $S_{1},...,S_{n}$ by the relations
\begin{eqnarray*}
\mathbf{P}_{x}^{+}(B_{n}) &:&=\frac{1}{V^{-}(x)}\mathbf{E}_{x}\left[
V^{-}(S_{n})I\left\{ B_{n}\right\} ;L_{n}\geq 0\right] \\
&=&\frac{1}{\hat{V}^{-}(x)}\mathbf{E}_{x}\left[ \hat{V}^{-}(S_{n})I\left\{
B_{n}\right\} ;L_{n}\geq 0\right]
\end{eqnarray*}%
and
\begin{eqnarray}
\mathbf{P}_{-x}^{-}(B_{n}) &:&=\frac{1}{V^{+}(x)}\mathbf{E}_{-x}\left[
V^{+}(-S_{n})I\left\{ B_{n}\right\} ;M_{n}<0\right]  \notag \\
&=&\frac{1}{\hat{V}^{+}(x)}\mathbf{E}_{-x}\left[ \hat{V}^{+}(-S_{n})I\left\{
B_{n}\right\} ;M_{n}<0\right] ,  \label{StayNegative}
\end{eqnarray}%
where $I\left\{ B_{n}\right\} $ is the indicator of the event $B_{n}$ and
the second equalities in the definitions are justified by (\ref%
{RenewRelation}). The expectations taken with respect to the measures $%
\mathbf{P}_{x}^{+}\left( \cdot \right) $ and $\mathbf{P}_{-x}^{-}\left(
\cdot \right) $ will be denoted by $\mathbf{E}_{x}^{+}\left[ \cdot \right] $
and $\mathbf{E}_{-x}^{-}\left[ \cdot \right] $, respectively.

In what follows to simplify the formulation of statements we use the
notation $k\bot r$, where $k\bot r=r$, if the parameters $r,k,$ and $n$ vary
in accordance with points 1) and 2) of Theorem \ref{T_beginningNew}, and $%
k\bot r=k$, if the parameters $r,k,$ and $n$ vary in accordance with points
3) or 4) of Theorem \ref{T_beginningNew}.

Let $\mathcal{F}_{l},l\in \mathbb{N}$, be the $\sigma $-algebra generated by
the random generating functions $F_{1}(s),F_{2}(s),...,F_{l}(s),$ and let $%
H_{1},H_{2},...,$ be a uniformly bounded sequence of random variables
adapted to the filtration $\mathcal{\tilde{F}=}\left\{ \mathcal{F}_{l},l\in
\mathbb{N}\right\} $ and converging $\mathbf{P}_{0}^{+}$-a.s. to a random
variable $H_{\infty }$ as $n\rightarrow \infty $. We take a parameter $%
k=k(n) $ growing to infinity as $n\rightarrow \infty $ in such a way that $%
k=o(n)$ and introduce, for $r<n$ and $m=n-r$ the conditional expectations%
\begin{equation*}
I_{n}^{\ast }(r,k;y,t):=\mathbf{E}\left[ H_{n}|S_{\tau _{r,n}}-S_{r}\leq
ya_{m},S_{n}\leq ta_{k},L_{n}\geq 0\right] ,y\in (-\infty ,0),
\end{equation*}%
and%
\begin{equation*}
I_{n}(r,k;y,t):=\mathbf{E}\left[ H_{n}|S_{\tau _{r,n}}\leq ya_{k\bot
r},S_{n}\leq ta_{k},L_{n}\geq 0\right] ,y\in \lbrack 0,\infty ).
\end{equation*}

The next theorem is a generalization of Lemma 2.5 in \cite{agkv} and Lemma 4
in \cite{VD2022}, where the definition
\begin{equation*}
\mathbf{P}_{x}^{+}(B)=\frac{1}{\hat{V}^{-}(x)}\mathbf{E}_{x}\left[ \hat{V}%
^{-}(S_{n})I\left\{ B\right\} ;L_{n}\geq 0\right]
\end{equation*}%
was used in the statements.

\begin{theorem}
\label{T_cond} Let Condition $B1$ be valid. Then

1) if $n\gg k\gg m=n-r$ then, for any $y\in (-\infty ,0)$ and $t>0$
\begin{equation}
\lim_{m\rightarrow \infty }I_{n}^{\ast }(r,k;y,t)=\mathbf{E}_{0}^{+}\left[
H_{\infty }\right] ;  \label{Cond3}
\end{equation}

2) if either $k\sim \theta m,\theta >0,$ or $\min (r,n-r)\gg k$, or $k\sim
\theta _{1}r,\theta _{1}>0$ or $k\gg r$ then, for any $y>0$ and $t>0$%
\begin{equation}
\lim_{k\wedge r\rightarrow \infty }I_{n}(r,k;y,t)=\mathbf{E}_{0}^{+}\left[
H_{\infty }\right] .  \label{Cond22}
\end{equation}
\end{theorem}

\textbf{Proof}. Our arguments mainly follow the scheme of proving Lemma 2.5
in \cite{agkv}. First we demonstrate the validity of (\ref{Cond3}). We fix $%
\lambda \geq 1,$ and agree to consider the expressions of the form $%
\lambda n,(\lambda -1)n$ and $\lambda r,\lambda m,$ as $\left[ \lambda n%
\right] ,\left[ (\lambda -1)n\right] $ and $\left[ \lambda r\right] ,\left[
\lambda m\right] .$ Introduce the notations
\begin{equation*}
\Psi \left( \lambda r,\lambda n,ya_{m},ta_{k}\right) :=\left\{ S_{\tau
_{\lambda r,\lambda n}}-S_{\lambda r}\leq ya_{m},S_{\lambda n}\leq
ta_{k},L_{n\lambda }\geq 0\right\}
\end{equation*}%
and%
\begin{equation*}
Q(y)=\mathbf{P}\left( \min_{0\leq s\leq 1}Y_{s}\leq y\right) .
\end{equation*}%
We show that, for any fixed $l\in \mathbb{N}$ and $y<0,$
\begin{equation}
\lim_{n\gg k\gg m\rightarrow \infty }\mathbf{E}\left[ H_{l}|\Psi \left(
r,n,ya_{m},ta_{k}\right) \right] =\mathbf{E}_{0}^{+}\left[ H_{l}\right] .
\label{Cond0}
\end{equation}%
To this aim we write the representation
\begin{eqnarray}
&&\mathbf{E}\left[ H_{l}|\Psi \left( r,n,ya_{m},ta_{k}\right) \right]  \notag
\\
&&\quad =\mathbf{E}\left[ H_{l}\frac{\mathbf{P}_{S_{l}}(S_{\tau
_{r-l,n-l}}^{\prime }-S_{r-l}^{\prime }\leq ya_{m},S_{n-l}^{\prime }\leq
ta_{k},L_{n-l}^{\prime }\geq 0)}{\mathbf{P}\left( \mathcal{B}%
(ta_{k},n)\right) };L_{l}\geq 0\right]  \notag \\
&&\qquad \qquad \times \frac{1}{\mathbf{P}\left( S_{\tau _{r,n}}-S_{r}\leq
ya_{m}|\mathcal{B}(ta_{k},n)\right) },  \label{Cond10}
\end{eqnarray}%
where $\mathcal{S}^{\prime }=\left\{ S_{n}^{\prime },n=0,1,2,...\right\} $
is a probabilistic copy of the random walk $\mathcal{S}$, being independent
of the random variables $\left\{ S_{j},j=0,1,...,l\right\} ,$ and let $%
L_{n}^{\prime }=\min {}_{0\leq i\leq n}$ $S_{i}^{\prime }.$ Note, that by
Theorem \ref{T_Minim} with $m=n-r$ we have
\begin{equation}
\lim_{n\gg k\gg m\rightarrow \infty }\mathbf{P}\left( S_{\tau
_{r,n}}-S_{r}\leq ya_{m}|\mathcal{B}(ta_{k},n)\right) =\mathbf{P}\left(
\min_{0\leq s\leq 1}Y_{s}\leq y\right) =Q(y).\   \label{Cond11}
\end{equation}%
Moreover, in view of (\ref{BasicAsymptotic})%
\begin{equation}
\lim_{n\rightarrow \infty }\frac{\mathbf{P}_{w}\left( \mathcal{B}%
(ta_{k},n-l)\right) }{\mathbf{P}\left( \mathcal{B}(ta_{k},n)\right) }%
=V^{-}\left( w\right)  \label{3.6af}
\end{equation}%
and, by (\ref{BasicAsymptotic}) and (\ref{Rough1}) there exist constants $%
C_{1},C_{2},$ and $C_{3}$ such that, for any fixed $l$ and all $n\geq l$
\begin{eqnarray}
&&\frac{\mathbf{P}_{w}(S_{\tau _{r-l,n-l}}^{\prime }-S_{r-l}^{\prime }\leq
ya_{m},S_{n-l}^{\prime }\leq ta_{k},L_{n-l}^{\prime }\geq 0)}{\mathbf{P}%
\left( \mathcal{B}(ta_{k},n)\right) }\leq \frac{\mathbf{P}_{w}(\mathcal{B}%
(ta_{k},n-l))}{\mathbf{P}\left( \mathcal{B}(ta_{k},n)\right) }  \notag \\
&&\qquad \qquad \qquad \qquad \quad \leq \frac{C_{1}\,b_{n-l}\,V^{-}(w)%
\int_{0}^{ta_{k}}\text{$V$}^{+}(z)dz}{C_{2}b_{n}\int_{0}^{ta_{k}}\text{$V$}%
^{+}(z)dz}\leq C_{3}V^{-}(w).  \label{Cond33}
\end{eqnarray}%
Further, applying (\ref{BasicAsymptotic}) and using definition (\ref{Defb})
and Theorem \ref{T_Minim}, we see, that, for each fixed $w\geq 0$ and $l\in
\mathbb{N}$
\begin{equation}
\lim_{n\gg k\gg m\rightarrow \infty }\frac{\mathbf{P}_{w}(S_{\tau
_{r-l,n-l}}^{\prime }-S_{r-l}^{\prime }\leq ya_{m},S_{n-l}^{\prime }\leq
ta_{k},L_{n-l}^{\prime }\geq 0)}{\mathbf{P}\left( \mathcal{B}%
(ta_{k},n)\right) }=Q(y)V^{-}(w).\text{ }  \label{Shift}
\end{equation}%
If the distribution of the random variable $X_{1}$ is absolutely continuous,
then $V^{-}(0)=1$ in view of (\ref{V_zero})\ and, therefore,
\begin{equation}
\mathbf{E}\left[ H_{l}V^{-}(S_{l});L_{l}\geq 0\right] =\frac{1}{V^{-}(0)}%
\mathbf{E}\left[ H_{l}V^{-}(S_{l});L_{l}\geq 0\right] \mathbf{=E}_{0}^{+}%
\left[ H_{l}\right] <\infty .  \label{Finite_expectation}
\end{equation}%
Hence, using the dominated convergence theorem, (\ref{3.6af}) and (\ref%
{Shift}), we conclude that, for each fixed $l\in \mathbb{N}$ and $y<0$
\begin{eqnarray*}
&&\lim_{n\gg k\gg m\rightarrow \infty }\mathbf{E}\left[ H_{l}\frac{\mathbf{P}%
_{S_{l}}(S_{\tau _{r-l,n-l}}^{\prime }-S_{r-l}^{\prime }\leq
ya_{m},S_{n-l}^{\prime }\leq ta_{k},L_{n-l}^{\prime }\geq 0)}{\mathbf{P}%
\left( \mathcal{B}(ta_{k},n)\right) };L_{l}\geq 0\right] \\
&&\quad =\mathbf{E}\left[ H_{l}\lim_{n\gg k\gg m\rightarrow \infty }\frac{%
\mathbf{P}_{S_{l}}(S_{\tau _{r-l,n-l}}^{\prime }-S_{r-l}^{\prime }\leq
ya_{m},S_{n-l}^{\prime }\leq ta_{k},L_{n-l}^{\prime }\geq 0)}{\mathbf{P}%
\left( \mathcal{B}(ta_{k},n)\right) };L_{l}\geq 0\right] \\
&&\quad =Q(y)\mathbf{E}\left[ H_{l}V^{-}(S_{l});L_{l}\geq 0\right] =Q(y)%
\mathbf{E}_{0}^{+}\left[ H_{l}\right] .
\end{eqnarray*}%
Combining this relation with (\ref{Cond10}) and (\ref{Cond11}), we get (\ref%
{Cond0}).

Now we prove that, for any $\lambda >1$
\begin{equation}
\lim_{n\gg k\gg m\rightarrow \infty }\mathbf{E}\left[ H_{n}|\Psi \left(
\lambda r,\lambda n,ya_{\lambda m},ta_{k}\right) \right] =\mathbf{E}_{0}^{+}%
\left[ H_{\infty }\right] .  \label{Cond44}
\end{equation}%
Since the relation $n\gg k\gg m\rightarrow \infty $ implies $\lambda n\gg
k\gg \lambda m\rightarrow \infty $, we conclude similarly to (\ref{Cond0})
that if $\lambda >0$ is fixed and $\lambda n\gg k\gg \lambda m\rightarrow
\infty $ then
\begin{equation}
\mathbf{E}\left[ H_{l}|\Psi \left( \lambda r,\lambda n,ya_{\lambda
m},ta_{k}\right) \right] \sim \mathbf{E}_{0}^{+}\left[ H_{l}\right] .
\label{vaf1}
\end{equation}

Further, in view of (\ref{Rough1}) the following estimates are valid
\begin{eqnarray*}
&&\left\vert \mathbf{E}\left[ \left( H_{n}-H_{l}\right) ;\Psi \left( \lambda
r,\lambda n,ya_{\lambda m},ta_{k}\right) \right] \right\vert \leq \mathbf{E}%
\left[ \left\vert H_{n}-H_{l}\right\vert ;S_{\lambda n}\leq
ta_{k},L_{n\lambda }\geq 0\right] \\
&&\qquad \qquad \qquad =\mathbf{E}\left[ \left\vert H_{n}-H_{l}\right\vert
\mathbf{P}_{S_{n}}(S_{(\lambda -1)n}^{\prime }\leq ta_{k},L_{n(\lambda
-1)}^{\prime }\geq 0);L_{n}\geq 0\right] \\
&&\qquad \qquad \qquad \leq Cb_{n(\lambda -1)}\int_{0}^{ta_{k}}\text{$V$}%
^{+}(z)dz\times \mathbf{E}\left[ \left\vert H_{n}-H_{l}\right\vert
\,\,V^{-}(S_{n});L_{n}\geq 0\right] \\
&&\qquad \qquad \qquad =Cb_{n(\lambda -1)}\int_{0}^{ta_{k}}\text{$V$}%
^{+}(z)dz\times \mathbf{E}_{0}^{+}\left[ \left\vert H_{n}-H_{l}\right\vert \,%
\right] .
\end{eqnarray*}%
Hence, using (\ref{Defb})--(\ref{Rough1}) with $w=0,$ we get
\begin{eqnarray*}
&&\frac{\left\vert \mathbf{E}\left[ (H_{n}-H_{l});\Psi \left( \lambda
r,\lambda n,ya_{\lambda m},ta_{k}\right) \right] \right\vert }{\mathbf{P}%
\left( \Psi \left( \lambda r,\lambda n,ya_{\lambda m},ta_{k}\right) \right) }
\\
&&\qquad \qquad \leq C\frac{b_{n(\lambda -1)}\int_{0}^{ta_{k}}\text{$V$}%
^{+}(z)dz}{C_{1}\,b_{n\lambda }\,\int_{0}^{ta_{k}}\text{$V$}^{+}(z)dz}%
\mathbf{E}_{0}^{+}\left[ \left\vert H_{n}-H_{l}\right\vert \right] \\
&&\qquad \qquad \leq C_{2}\left( \frac{\lambda }{\lambda -1}\right)
^{1+1/\alpha }\mathbf{E}_{0}^{+}\left[ \left\vert H_{n}-H_{l}\right\vert %
\right] .
\end{eqnarray*}%
Since the elements of the sequence $\{H_{k},k=1,2,...\}$ are uniformly
bounded we may apply the dominated convergence theorem to conclude that
\begin{equation*}
\lim_{l\rightarrow \infty }\lim_{n\rightarrow \infty }\mathbf{E}_{0}^{+}%
\left[ \left\vert H_{n}-H_{l}\right\vert \right] =0.
\end{equation*}%
Hence, using (\ref{vaf1}), we get, for any $\lambda >1$
\begin{eqnarray*}
&&\lim_{n\gg k\gg m\rightarrow \infty }\mathbf{E}\left[ H_{n}|\Psi \left(
\lambda r,\lambda n,ya_{\lambda m},ta_{k}\right) \right] \\
&&\qquad =\lim_{l\rightarrow \infty }\lim_{n\gg k\gg m\rightarrow \infty }%
\frac{\mathbf{E}\left[ \left( H_{n}-H_{l}\right) ;\Psi \left( \lambda
r,\lambda n,ya_{\lambda m},ta_{k}\right) \right] }{\mathbf{P}\left( \Psi
\left( \lambda r,\lambda n,ya_{\lambda m},ta_{k}\right) \right) } \\
&&\qquad \quad +\lim_{l\rightarrow \infty }\lim_{n\gg k\gg m\rightarrow
\infty }\frac{\mathbf{E}\left[ H_{l};\Psi \left( \lambda r,\lambda
n,ya_{\lambda m},ta_{k}\right) \right] }{\mathbf{P}\left( \Psi \left(
\lambda r,\lambda n,ya_{\lambda m},ta_{k}\right) \right) } \\
&&\qquad =\lim_{l\rightarrow \infty }\mathbf{E}_{0}^{+}\left[ H_{l}\right] =%
\mathbf{E}_{0}^{+}\left[ H_{\infty }\right] ,
\end{eqnarray*}%
which proves (\ref{Cond44}). To prove (\ref{Cond3}) it remains to show that
we may take $\lambda =1$ in the estimate given above. To justify such a
transition it is convenient to rewrite equality (\ref{Cond44}) in the form
\begin{eqnarray*}
&&\mathbf{E}\left[ H_{n};\Psi \left( \lambda r,\lambda n,ya_{\lambda
m},ta_{k}\right) \right] \\
&&\qquad =\left( \mathbf{E}_{0}^{+}\left[ H_{\infty }\right] +o(1)\right)
\mathbf{P}\left( \Psi \left( \lambda r,\lambda n,ya_{\lambda
m},ta_{k}\right) \right)
\end{eqnarray*}%
as $n\gg k\gg m\rightarrow \infty $. Assuming without loss of generality
that $H_{\infty }>0$ $\mathbf{P}_{0}^{+}$-a.s.\ and \ $\mathbf{E}_{0}^{+}%
\left[ H_{\infty }\right] \leq 1,$ we conclude that
\begin{eqnarray}
&&|\mathbf{E}\left[ H_{n};\Psi \left( r,n,ya_{m},ta_{k}\right) \right] -%
\mathbf{E}_{0}^{+}\left[ H_{\infty }\right] \mathbf{P}\left( \Psi \left(
\lambda r,\lambda n,ya_{\lambda m},ta_{k}\right) \right) |  \notag \\
&&\quad \leq |\mathbf{E}\left[ H_{n};\Psi \left( \lambda r,\lambda
n,ya_{\lambda m},ta_{k}\right) \right] -\mathbf{E}_{0}^{+}\left[ H_{\infty }%
\right] \mathbf{P}\left( \Psi \left( \lambda r,\lambda n,ya_{\lambda
m},ta_{k}\right) \right) |  \notag \\
&&\qquad +|\mathbf{P}\left( \Psi \left( \lambda r,\lambda n,ya_{\lambda
m},ta_{k}\right) \right) -\mathbf{P}\left( \Psi \left(
r,n,ya_{m},ta_{k}\right) \right) |.  \label{Cond100}
\end{eqnarray}%
By (\ref{Cond44}) the first difference on the right-hand part of this
inequality is
\begin{equation}
o\left( \mathbf{P}\left( \Psi \left( \lambda r,\lambda n,ya_{\lambda
m},ta_{k}\right) \right) \right)  \label{Cond110}
\end{equation}%
as $n\gg k\gg m\rightarrow \infty $, and, therefore, is $o\left( \mathbf{P}%
\left( S_{n}\leq ta_{k},L_{n}\geq 0\right) \right) $, since
\begin{equation}
\lim_{n\gg k\rightarrow \infty }\frac{\mathbf{P}\left( S_{n\lambda }\leq
ta_{k},L_{n\lambda }\geq 0\right) }{\mathbf{P}\left( S_{n}\leq
ta_{k},L_{n}\geq 0\right) }=\lim_{n\rightarrow \infty }\frac{b_{n\lambda }}{%
b_{n}}=\lambda ^{1+1/\alpha }  \label{RatioB}
\end{equation}%
by (\ref{BasicAsymptotic}) and (\ref{Defb}).

Further, using (\ref{BasicAsymptotic}) once more and Theorem \ref{T_Minim}
and definition (\ref{Defb}), we obtain
\begin{eqnarray*}
&&|\mathbf{P}\left( \Psi \left( \lambda r,\lambda n,ya_{\lambda
m},ta_{k}\right) \right) -\mathbf{P}\left( \Psi \left(
r,n,ya_{m},ta_{k}\right) \right) | \\
&\leq &\left\vert \mathbf{P}\left( \Psi \left( \lambda r,\lambda
n,ya_{\lambda m},ta_{k}\right) \right) -g_{\alpha ,\beta }(0)b_{n\lambda
}\int_{0}^{ta_{k}}\text{$V$}^{+}(z)dz\right\vert \\
&&\quad +\left\vert \mathbf{P}\left( \Psi \left( r,n,ya_{m},ta_{k}\right)
\right) -g_{\alpha ,\beta }(0)b_{n}\int_{0}^{ta_{k}}\text{$V$}%
^{+}(z)dz\right\vert \\
&&\quad +g_{\alpha ,\beta }(0)\left\vert b_{n\lambda }-b_{n}\right\vert
\int_{0}^{ta_{k}}\text{$V$}^{+}(z)dz \\
&=&o\left( b_{n}\int_{0}^{ta_{k}}\text{$V$}^{+}(z)dz\right) +g_{\alpha
,\beta }(0)b_{n}\left\vert \frac{b_{n\lambda }}{b_{n}}-1\right\vert
\int_{0}^{ta_{k}}\text{$V$}^{+}(z)dz.
\end{eqnarray*}%
Hence, passing to the limit as $n\gg k\gg m\rightarrow \infty ,$ letting $%
\lambda $ to~$1$ and using (\ref{RatioB}), we deduce that
\begin{equation}
\lim_{\lambda \downarrow 1}\limsup_{n\gg k\gg m\rightarrow \infty }\frac{%
\left\vert \mathbf{P}\left( \Psi \left( \lambda r,\lambda
n,ya_{m},ta_{k}\right) \right) -\mathbf{P}\left( \Psi \left(
r,n,ya_{m},ta_{k}\right) \right) \right\vert }{b_{n}\int_{0}^{ta_{k}}\text{$%
V $}^{+}(z)dz}=0.  \label{Cond101}
\end{equation}%
Combining (\ref{Cond100})--(\ref{Cond101}) with (\ref{Cond10}) we get (\ref%
{Cond3}).

We now prove (\ref{Cond22}). For fixed $1\leq l<n$ and $y\in (0,\infty )$ we
consider the conditional expectation
\begin{eqnarray*}
&&\mathbf{E}\left[ H_{l}|S_{\tau _{r,n}}\leq ya_{k\bot r},S_{n}\leq
ta_{k},L_{n}\geq 0\right] \\
&&\quad =\mathbf{E}\left[ H_{l}\frac{\mathbf{P}_{S_{l}}(S_{\tau
_{r-l,n-l}}^{\prime }\leq ya_{k\bot r},S_{n-l}^{\prime }\leq
ta_{k},L_{n-l}^{\prime }\geq 0)}{\mathbf{P}\left( \mathcal{B}%
(ta_{k},n)\right) };L_{l}\geq 0\right] \\
&&\times \frac{1}{\mathbf{P}\left( S_{\tau _{r,n}}\leq ya_{k\bot r}|\mathcal{%
B}(ta_{k},n)\right) }.
\end{eqnarray*}%
We know by (\ref{BasicAsymptotic}) and Theorem \ref{T_beginningNew} that
\begin{equation*}
\mathbf{P}\left( S_{\tau _{r,n}}\leq ya_{k\bot r}|\mathcal{B}%
(ta_{k},n)\right) \sim \hat{Q}(y),
\end{equation*}%
where
\begin{equation}
\hat{Q}(y)=\left\{
\begin{array}{ccc}
C^{\ast \ast }\mathcal{H}(y) & \text{if} & n\gg k\gg r\rightarrow \infty ,
\\
&  &  \\
W(\theta ^{1/\alpha }t,\theta ^{1/\alpha }t\wedge y) & \text{if} & n\gg
k\sim \theta r\rightarrow \infty , \\
&  &  \\
1-\left( 1-\frac{t\wedge y}{t}\right) ^{\alpha \rho +1} & \text{if } & \min
(r,n-r)\gg k\rightarrow \infty , \\
&  &  \\
A\left( \theta ^{1/\alpha }t,\theta ^{1/\alpha }\left( t\wedge y\right)
\right) & \text{if} & n\gg k\sim \theta (n-r)\rightarrow \infty .%
\end{array}%
\right.  \label{DefQ^}
\end{equation}

Further, using the same arguments we have applied to prove (\ref{Cond33}),
we conclude that there exists a constant $C_{2}$ such, that, for any fixed $%
l\in \mathbb{N}$, all $n\geq l$ and $y\in (0,t]$
\begin{equation*}
\frac{\mathbf{P}_{w}(S_{\tau _{r-l,n-l}}^{\prime }\leq ya_{k\bot
r},S_{n-l}^{\prime }\leq ta_{k},L_{n-l}^{\prime }\geq 0)}{\mathbf{P}\left(
S_{\tau _{r,n}}\leq ya_{k\bot r}|\mathcal{B}(ta_{k},n)\right) }\leq
C_{2}V^{-}(w).
\end{equation*}%
Using (\ref{BasicAsymptotic}), definition (\ref{Defb}) and applying Theorem~%
\ref{T_beginningNew}, we see that, for any fixed $w\geq 0,t>0$ and $l\in
\mathbb{N}$
\begin{equation}
\frac{\mathbf{P}_{w}(S_{\tau _{r-l,n-l}}^{\prime }\leq ya_{k\bot
r},S_{n-l}^{\prime }\leq ta_{k},L_{n-l}^{\prime }\geq 0)}{\mathbf{P}\left(
S_{\tau _{r,n}}\leq ya_{k\bot r}|\mathcal{B}(ta_{k},n)\right) }\sim V^{-}(w)%
\hat{Q}(y)\text{ }  \label{Shift2}
\end{equation}%
if $n\rightarrow \infty $ and the parameters $r$\ and $k$ vary within the
ranges specified by (\ref{DefQ^}). This estimate combined with (\ref%
{Finite_expectation}) allows us to apply the dominated convergence theorem
and conclude in view of (\ref{Shift2}) that, for any fixed $l\in \mathbb{N}$
and $y\geq 0$
\begin{eqnarray*}
&&\lim_{n\rightarrow \infty }\mathbf{E}\left[ H_{l}\frac{\mathbf{P}%
_{S_{l}}(S_{\tau _{r-l,n-l}}^{\prime }\leq ya_{k\bot r},S_{n-l}^{\prime
}\leq ta_{k},L_{n-l}^{\prime }\geq 0)}{\mathbf{P}\left( S_{n}\leq
ta_{k},L_{n}\geq 0\right) };L_{l}\geq 0\right] \\
&&\quad =\mathbf{E}\left[ H_{l}\lim_{n\rightarrow \infty }\frac{\mathbf{P}%
_{S_{l}}(S_{\tau _{r-l,n-l}}^{\prime }\leq ya_{k\bot r},S_{n-l}^{\prime
}\leq ta_{k},L_{n-l}^{\prime }\geq 0)}{\mathbf{P}\left( S_{n}\leq
ta_{k},L_{n}\geq 0\right) };L_{l}\geq 0\right] \\
&&\quad =\frac{1}{V^{-}(0)}\mathbf{E}\left[ H_{l}V^{-}(S_{l});L_{l}\geq 0%
\right] \hat{Q}(y)=\mathbf{E}_{0}^{+}\left[ H_{l}\right] \hat{Q}(y).
\end{eqnarray*}

To complete the proof of (\ref{Cond22}) it is now necessary to introduce a
temporary notation
\begin{equation*}
\hat{\Psi}\left( \lambda r,\lambda n,ya_{k\bot r},ta_{k}\right) :=\left\{
S_{\tau _{\lambda r,\lambda n}}\leq ya_{k\bot r},S_{\lambda n}\leq
ta_{k},L_{n\lambda }\geq 0\right\} ,
\end{equation*}%
for fixed $\lambda \geq 1$ and then to repeat almost literary the steps used
at the end of proving (\ref{Cond3}) starting from formula (\ref{Cond44}) and
replacing, where it is necessary, $Q$ by $\hat{Q}$.

Theorem \ref{T_cond} is proved.

The next statement easily follows from Theorem \ref{T_cond} and does not
require the boundness of the random variables involved.

\begin{corollary}
\label{C_cond} Assume that Condition $B1$ is valid. Let $H_{1},H_{2},...,$
be a nonnegative sequence of random variables adapted to the filtration $%
\mathcal{\tilde{F}=}\left\{ \mathcal{F}_{l},l\in \mathbb{N}\right\} ,$ and
converging $\mathbf{P}_{0}^{+}$-a.s. to a random variable $H_{\infty }$ as $%
n\rightarrow \infty $. If either $k\sim \theta m,\theta >0,$ or $\min
(r,n-r)\gg k$, or $k\sim \theta r,\theta >0,$ or $k\gg r,$ then, for all
positive $y,t,$ and $\lambda $%
\begin{equation*}
\lim_{k\bot r\rightarrow \infty }\mathbf{E}\left[ e^{-\lambda H_{n}}|S_{\tau
_{r,n}}\leq ya_{k\bot r},S_{n}\leq ta_{k},L_{n}\geq 0\right] =\mathbf{E}%
_{0}^{+}\left[ e^{-\lambda H_{\infty }}\right] ,
\end{equation*}%
and, therefore,
\begin{equation*}
\lim_{k\bot r\rightarrow \infty }\mathbf{P}\left( H_{n}>x|S_{\tau
_{r,n}}\leq ya_{k\bot r},S_{n}\leq ta_{k},L_{n}\geq 0\right) =\mathbf{P}%
_{0}^{+}\left( H_{\infty }>x\right)
\end{equation*}%
for any continuity point of the distribution of $H_{\infty }$.
\end{corollary}

We note that the statement of the corollary still holds if we use either of
the following conditions: $\left\{ S_{n}\leq ta_{k},L_{n}\geq 0\right\} $ or
$\left\{ L_{n}\geq 0\right\} ,$ instead of the original condition $\left\{
S_{\tau _{r,n}}\leq ya_{k\bot r},S_{n}\leq ta_{k},L_{n}\geq 0\right\} $ (see
Lemma 4 in \cite{VD2022} and Lemma 2.5 in \cite{agkv}, respectively).
Observe also that under the condition $\left\{ S_{n}\leq ta_{k},L_{n}\geq
0\right\} $ it is necessary to pass to the limit as $n\gg k\rightarrow
\infty $, and in the case of the condition $\left\{ L_{n}\geq 0\right\} $ --
to the limit as $n\rightarrow \infty $.

\section{Change of condition}

In this section we replace the condition $\mathcal{B}\left( ta_{k},n\right)
=\left\{ S_{n}\leq ta_{k},L_{n}\geq 0\right\} ,$ used in Theorems \ref%
{T_beginningNew} and \ref{T_Minim} and does not \ directly connected with
the properties of branching processes in random environment, by the
condition $\mathcal{R}\left( ta_{k},n\right) =\left\{ S_{n}\leq
ta_{k},Z_{n}>0\right\} $. Such a replacement allows us to prove in the
section \ref{Sec5} limit Theorems \ref{T_reducedSmall}-\ref{T_GeneralFrac},
describing asymptotic properties of reduced critical branching processes in
random environment.

Set
\begin{equation*}
F_{l,n}(s):=F_{l+1}(F_{l+2}(\ldots (F_{n}(s))\ldots )),\;0\leq l\leq
n-1,\;F_{n,n}(s):=s,
\end{equation*}%
and denote
\begin{equation*}
F_{l,\infty }(s):=\lim_{n\rightarrow \infty }F_{l,n}(s),0\leq s\leq 1.
\end{equation*}%
Note, that the last limit always exists (see Theorem 5 in \cite{AK1971}).

Denote $\mathcal{A}_{u.s.}$ the event $\{Z_{n}>0$ for all $n>0\}.$ Our
starting point is the following relation established in Theorem 1 in \cite%
{VD2022} under the conditions B1 and B2: if $n\gg k\rightarrow \infty ,$
then, for any $t>0$
\begin{eqnarray}
\mathbf{P}\left( S_{n}\leq ta_{k},Z_{n}>0\right) &\sim &\Theta \mathbf{P}%
\left( S_{n}\leq ta_{k},L_{n}\geq 0\right)  \notag \\
&\sim &\Theta g_{\alpha ,\beta }(0)b_{n}\int_{0}^{ta_{k}}V^{+}(u)du,
\label{AsympMain1}
\end{eqnarray}%
where
\begin{eqnarray}
\Theta &=&\sum_{j=0}^{\infty }\sum_{i=1}^{\infty }\mathbf{P}(Z_{j}=i,\tau
_{j}=j)\mathbf{E}_{0}^{+}\left[ 1-F_{0,\infty }^{i}(0)\right]  \notag \\
&=&\sum_{j=0}^{\infty }\sum_{i=1}^{\infty }\mathbf{P}(Z_{j}=i,\tau _{j}=j)%
\mathbf{P}_{0}^{+}\left( \mathcal{A}_{u.s.}|Z_{0}=i\right) .
\label{DefTheta}
\end{eqnarray}

Note that
\begin{eqnarray}
0 &<&\Theta \leq \sum_{j=0}^{\infty }\sum_{i=1}^{\infty }\mathbf{P}%
(Z_{j}=i,\tau _{j}=j)=\sum_{j=0}^{\infty }\mathbf{P}\left( Z_{j}>0,\tau
_{j}=j\right)  \notag \\
&\leq &\sum_{j=0}^{\infty }\mathbf{E}\left[ e^{S_{j}},\tau _{j}=j\right]
=\sum_{j=0}^{\infty }\mathbf{E}\left[ e^{S_{j}},M_{j}<0\right] <\infty ,
\label{Sparr}
\end{eqnarray}%
since
\begin{equation}
\mathbf{E}\left[ e^{S_{j}},M_{j}<0\right] \sim Cb_{j}=\frac{C}{j^{1/\alpha
+1}l(j)},\quad C\in (0,\infty ),  \label{ExpFunct}
\end{equation}%
as $j\rightarrow \infty $ by Proposition 2.1 in \cite{ABGV2011} and
definition (\ref{Defb}).

\begin{lemma}
\label{L_preliminary} Assume that conditions $B1$ and $B2$ are valid.$ $Let
$\mathcal{T}_{n}$ be an event from the $\sigma $-algebra generated by the
set $F_{1},...,F_{n},Z_{0},Z_{1},...,Z_{n}$. Then, for any \ fixed $w\geq 0$
and $t>0$

1)
\begin{equation}
\lim_{J\rightarrow \infty }\limsup_{n\gg k\rightarrow \infty }\frac{\mathbf{P%
}_{w}\left( \mathcal{T}_{n},S_{n}\leq ta_{k},Z_{n}>0,\tau _{n}>J\right) }{%
\mathbf{P}\left( S_{n}\leq ta_{k},L_{n}\geq 0\right) }=0;  \label{Prelim1}
\end{equation}%
2) for any fixed $J\in \mathbb{N}$%
\begin{equation}
\limsup_{n\gg k\rightarrow \infty }\frac{\mathbf{P}_{w}\left( \mathcal{T}%
_{n},S_{n}\leq ta_{k},Z_{n}>0,\tau _{n}\leq J,S_{\tau _{n}}\leq -\sqrt{ta_{k}%
}\right) }{\mathbf{P}\left( S_{n}\leq ta_{k},L_{n}\geq 0\right) }=0
\label{Prelim2}
\end{equation}%
and for any fixed $j\in \lbrack 0,J]$
\begin{equation}
\lim_{K\rightarrow \infty }\limsup_{n\gg k\rightarrow \infty }\frac{\mathbf{P%
}_{w}\left( \mathcal{T}_{n},S_{n}\leq ta_{k},Z_{j}>K,Z_{n}>0,\tau
_{n}=j,S_{j}>-\sqrt{ta_{k}}\right) }{\mathbf{P}\left( S_{n}\leq
ta_{k},L_{n}\geq 0\right) }=0.  \label{Prelim3}
\end{equation}%
(Here $S_{n}-S_{n-1}=\log F_{n}^{^{\prime }}\left( 1\right) ,n\geq 1$).
\end{lemma}

\textbf{Proof}. Using Lemma 5 in \cite{VD2022} and~(\ref{BasicAsymptotic})
we conclude that
\begin{eqnarray*}
&&\lim_{J\rightarrow \infty }\limsup_{n\gg k\rightarrow \infty }\frac{%
\mathbf{P}_{w}\left( S_{n}\leq ta_{k},Z_{n}>0,\tau _{n}>J\right) }{\mathbf{P}%
\left( S_{n}\leq ta_{k},L_{n}\geq 0\right) } \\
&&\qquad =\lim_{J\rightarrow \infty }\limsup_{n\gg k\rightarrow \infty }%
\frac{\mathbf{P}\left( S_{n}\leq ta_{k}-w,Z_{n}>0,\tau _{n}>J\right) }{%
\mathbf{P}\left( S_{n}\leq ta_{k},L_{n}\geq 0\right) }=0.
\end{eqnarray*}%
This relation implies (\ref{Prelim1}).

Further, we have for any fixed $j\in \lbrack 0,J]$
\begin{eqnarray*}
&&\mathbf{P}_{w}\left( \mathcal{T}_{n},S_{n}\leq ta_{k},Z_{n}>0,\tau
_{n}=j,S_{j}\leq -\sqrt{ta_{k}}\right) \\
&&\qquad \leq \mathbf{P}\left( S_{n}\leq ta_{k}-w,Z_{n}>0,\tau
_{n}=j,S_{j}\leq -\sqrt{ta_{k}}-w\right) \\
&&\qquad \leq \mathbf{E}\left[ e^{S_{\tau _{n}}};S_{n}\leq ta_{k},\tau
_{n}=j,S_{j}\leq -\sqrt{ta_{k}}\right]
\end{eqnarray*}%
and (compare with formulas (7) and (8) in Lemma 5 in \cite{VD2022})
\begin{eqnarray*}
&&\mathbf{E}\left[ e^{S_{\tau _{n}}};S_{n}\leq ta_{k},\tau _{n}=j,S_{j}\leq -%
\sqrt{ta_{k}}\right] \\
&&\quad \leq Cb_{j}b_{n}\sum_{i\geq \sqrt{ta_{k}}}e^{-i}V^{-}(i)%
\int_{0}^{ta_{k}+i}V^{+}(u)du \\
&&\quad \leq Cb_{j}b_{n}\sum_{i\geq \sqrt{ta_{k}}}e^{-i}V^{-}(i)\left(
\int_{0}^{2ta_{k}}V^{+}(u)du+\int_{0}^{2i}V^{+}(u)du\right) \\
&&\qquad =o\left( b_{n}\int_{0}^{ta_{k}}V^{+}(u)du\right)
\end{eqnarray*}%
as $n\rightarrow \infty ,$ since the series
\begin{equation*}
\sum_{i=1}^{\infty }e^{-i}V^{-}(i)\left( 1+\int_{0}^{2i}V^{+}(u)du\right)
\end{equation*}%
is convergent in view of (\ref{Regular1})--(\ref{Regular2}) and by (\ref%
{AsympV}) there exists a constant $C\in (0,\infty )$ such, that
\begin{equation*}
\int_{0}^{2ta_{k}}V^{+}(u)du\leq C\int_{0}^{ta_{k}}V^{+}(u)du
\end{equation*}%
for all $ta_{k}\geq 1$. Combining the obtained estimates proves (\ref%
{Prelim2}).

Finally, setting $S_{i}^{\prime }=S_{j+i}-S_{j},i=0,1,2,...$ and $%
L_{m}^{\prime }=\min_{0\leq i\leq m}S_{i}^{\prime },$ we have for any fixed $%
K\geq 1$
\begin{eqnarray*}
&&\mathbf{P}_{w}\left( \mathcal{T}_{n},S_{n}\leq ta_{k},Z_{n}>0,Z_{j}>K,\tau
_{n}=j,S_{j}>-\sqrt{ta_{k}}\right) \\
&&\qquad \leq \mathbf{P}\left( S_{n}\leq ta_{k}-w,\tau
_{n}=j,Z_{j}>K,S_{j}\geq -\sqrt{ta_{k}}-w\right) \\
&&\qquad \leq \mathbf{P}\left( S_{n}-S_{j}\leq ta_{k}+\sqrt{ta_{k}},\tau
_{n}=j,Z_{j}>K\right) \\
&&\qquad =\mathbf{P}\left( S_{n-j}^{\prime }\leq ta_{k}+\sqrt{ta_{k}}%
,L_{n-j}^{\prime }\geq 0,\tau _{j}=j,Z_{j}>K\right) \\
&&\qquad =\mathbf{P}\left( \tau _{j}=j,Z_{j}>K\right) \mathbf{P}\left(
S_{n-j}\leq ta_{k}+\sqrt{ta_{k}},L_{n-j}\geq 0\right) \\
&&\qquad \leq K^{-1}\mathbf{E}\left[ e^{S_{j}},\tau _{j}=j\right] \mathbf{P}%
\left( S_{n-j}\leq ta_{k}+\sqrt{ta_{k}},L_{n-j}\geq 0\right) \\
&&\qquad \leq K^{-1}\mathbf{P}\left( S_{n-j}\leq ta_{k}+\sqrt{ta_{k}}%
,L_{n-j}\geq 0\right) .
\end{eqnarray*}%
Hence, letting $K$ to infinity and using (\ref{BasicAsymptotic}), it is not
difficult to check the validity of (\ref{Prelim3}).

Lemma \ref{L_preliminary} is proved.

We apply Lemma \ref{L_preliminary} to establish below a number of important
statements where, as before, we use the notation $m=n-r$.

\begin{lemma}
\label{L_SmallSurvival} Assume that Conditions $B1$ and $B2$ are valid.
Then, for any $t>0,w>0$ and $y\leq 0$
\begin{equation*}
\lim_{n\gg k\gg m\rightarrow \infty }\mathbf{P}_{w}\left( \frac{S_{\tau
_{r,n}}-S_{r}}{a_{m}}\leq y\Big|S_{n}\leq ta_{k},Z_{n}>0\right) =\mathbf{P}%
\left( \min_{0\leq s\leq 1}Y_{s}\leq y\right) .
\end{equation*}
\end{lemma}

\textbf{Proof}. Note that in view of (\ref{AsympV}), for any $t>0$
\begin{equation}
\int_{0}^{ta_{k}}V^{+}(z)dz\sim \frac{1}{\alpha \rho +1}ta_{k}V^{+}(ta_{k})
\label{vaf2}
\end{equation}%
as $k\rightarrow \infty $. It follows from (\ref{AsympMain1})\ that, for
each fixed $w\geq 0$ and $t>0$
\begin{eqnarray}
&&\mathbf{P}_{w}\left( S_{n}\leq ta_{k},Z_{n}>0\right) =\mathbf{P}\left(
S_{n}\leq ta_{k}-w,Z_{n}>0\right)  \notag \\
&&\quad \sim \Theta g_{\alpha ,\beta }(0)b_{n}\int_{0}^{ta_{k}-w}\text{$V$}%
^{+}(z)dz\sim \Theta g_{\alpha ,\beta }(0)b_{n}\int_{0}^{ta_{k}}\text{$V$}%
^{+}(z)dz  \notag \\
&&\quad \sim \Theta \mathbf{P}\left( S_{n}\leq ta_{k},L_{n}\geq 0\right)
\label{Equival}
\end{eqnarray}%
as $n\gg k\rightarrow \infty $.\

We now fixed positive integers $K$ and $J<n$ and setting $\mathcal{T}%
_{r,n}(x)=\left\{ S_{\tau _{r,n}}-S_{r}\leq x\right\} $ write the
decomposition
\begin{eqnarray*}
&&\mathbf{P}_{w}\left( \mathcal{T}_{r,n}(ya_{m}),S_{n}\leq
ta_{k},Z_{n}>0\right) \\
&&\qquad =\mathbf{P}_{w}\left( \mathcal{T}_{r,n}(ya_{m}),S_{n}\leq
ta_{k},Z_{n}>0,\tau _{n}>J\right) \\
&&\qquad +\mathbf{P}_{w}\left( \mathcal{T}_{r,n}(ya_{m}),S_{n}\leq
ta_{k},Z_{n}>0,\tau _{n}\leq J,S_{\tau _{n}}\leq -\sqrt{ta_{k}}\right) \\
&&\qquad +\mathbf{P}_{w}\left( \mathcal{T}_{r,n}(ya_{m}),S_{n}\leq
ta_{k},Z_{n}>0,Z_{\tau _{n}}>K,\tau _{n}\leq J,S_{\tau _{n}}>-\sqrt{ta_{k}}%
\right) \\
&&\qquad +\mathbf{P}_{w}\left( \mathcal{T}_{r,n}(ya_{m}),S_{n}\leq
ta_{k},Z_{n}>0,Z_{\tau _{n}}\leq K,\tau _{n}\leq J,S_{\tau _{n}}>-\sqrt{%
ta_{k}}\right) .
\end{eqnarray*}

By Lemma \ref{L_preliminary} and relation (\ref{Equival}) we conclude that
to prove the lemma it is sufficient to show that
\begin{eqnarray}
&&\lim_{J\rightarrow \infty }\lim_{K\rightarrow \infty }\lim_{n\gg
k\rightarrow \infty }\frac{\mathbf{P}_{w}\left( \mathcal{T}%
_{r,n}(ya_{m}),S_{n}\leq ta_{k},Z_{n}>0,Z_{\tau _{n}}\leq K,\tau _{n}\leq
J,S_{\tau _{n}}>-\sqrt{ta_{k}}\right) }{\Theta \mathbf{P}\left( S_{n}\leq
ta_{k},L_{n}\geq 0\right) }  \notag \\
&&\qquad \qquad \qquad \qquad =\mathbf{P}\left( \min_{0\leq s\leq
1}Y_{s}\leq y\right) .  \label{Afan0}
\end{eqnarray}%
To check the validity of this statement we fix $j\in \lbrack 0,J],J<n,$ and
investigate the asymptotic behavior of the probability
\begin{eqnarray*}
&&\mathbf{P}_{w}\left( \mathcal{T}_{r,n}(ya_{m}),S_{n}\leq
ta_{k},Z_{n}>0,Z_{j}\leq K,\tau _{n}=j,S_{j}>-\sqrt{ta_{k}}\right) \\
&&\qquad =\int_{-\sqrt{ta}}^{w}\sum_{i=1}^{K}\mathbf{P}_{w}\left( S_{j}\in
dq,Z_{j}=i,\tau _{j}=j\right) \\
&&\times \mathbf{E}\left[ \mathbf{P}\left( Z_{n-j}>0|\mathcal{E}%
,Z_{0}=i\right) ;\mathcal{T}_{r-j,n-j}(ya_{m}),S_{n-j}\leq
ta_{k}-q,L_{n-j}\geq 0\right]
\end{eqnarray*}%
as $n\gg k\rightarrow \infty $. By the $sigma$-additivity axiom and monotonicity of the survival probability for each $i\in \mathbb{N}$ we have
\begin{eqnarray}
H_{n-j}(i):= &&\mathbf{P}\left( Z_{n-j}>0|\mathcal{E},Z_{0}=i\right)  \notag
\\
&\rightarrow &\mathbf{P}_{0}^{+}\left( \mathcal{A}_{u.s}|\mathcal{E}%
,Z_{0}=i\right) =:H_{\infty }\left( i\right)  \label{ConvAs}
\end{eqnarray}%
$\mathbf{P}_{0}^{+}$-a.s. as $n-j\rightarrow \infty $. Moreover, $H_{\infty
}\left( i\right) >0$ $\mathbf{P}_{0}^{+}$-a.s. according Proposition 3.1 in
\cite{agkv}. Further, for $q\in (-\sqrt{ta_{k}},w]$ we obtain
\begin{eqnarray*}
&&\mathbf{E}\left[ H_{n-j}(i);\mathcal{T}_{r-j,n-j}(ya_{m}),S_{n-j}\leq
ta_{k}+\sqrt{ta_{k}},L_{n-j}\geq 0\right] \\
&&\quad \geq \mathbf{E}\left[ H_{n-j}(i);\mathcal{T}%
_{r-j,n-j}(ya_{m}),S_{n-j}\leq ta_{k}-q,L_{n-j}\geq 0\right] \\
&&\quad \geq \mathbf{E}\left[ H_{n-j}(i);\mathcal{T}%
_{r-j,n-j}(ya_{m}),S_{n-j}\leq ta_{k}-w,L_{n-j}\geq 0\right] .
\end{eqnarray*}%
Note that in view of Theorem \ref{T_cond}
\begin{eqnarray*}
&&\mathbf{E}\left[ H_{n-j}(i);\mathcal{T}_{r-j,n-j}(ya_{m}),S_{n-j}\leq
ta_{k}-w,L_{n-j}\geq 0\right] \\
&&\qquad \sim \mathbf{P}\left( \mathcal{T}_{r-j,n-j}(ya_{m}),S_{n-j}\leq
ta_{k}-w,L_{n-j}\geq 0\right) \mathbf{E}_{0}^{+}\left[ H_{\infty }(i)\right]
\\
&&\qquad =\mathbf{P}\left( \mathcal{T}_{r-j,n-j}(ya_{m})|S_{n-j}\leq
ta_{k}-w,L_{n-j}\geq 0\right) \times \\
&&\qquad \quad \times \mathbf{P}\left( S_{n-j}\leq ta_{k}-w,L_{n-j}\geq
0\right) \mathbf{P}_{0}^{+}\left( \mathcal{A}_{u.s}|Z_{0}=i\right)
\end{eqnarray*}%
as $n\gg k\rightarrow \infty $. The similar arguments show that
\begin{eqnarray*}
&&\mathbf{E}\left[ H_{n-j}(i);\mathcal{T}_{r-j,n-j}(ya_{m}),S_{n-j}\leq
ta_{k}+\sqrt{ta_{k}},L_{n-j}\geq 0\right] \\
&&\quad \sim \mathbf{P}\left( \mathcal{T}_{r-j,n-j}(ya_{m}),S_{n-j}\leq
ta_{k}+\sqrt{ta_{k}},L_{n-j}\geq 0\right) \mathbf{P}_{0}^{+}\left( \mathcal{A%
}_{u.s}|Z_{0}=i\right) \\
&&\quad =\mathbf{P}\left( \mathcal{T}_{r-j,n-j}(ya_{m})|S_{n-j}\leq ta_{k}+%
\sqrt{ta_{k}},L_{n-j}\geq 0\right) \times \\
&&\quad \quad \times \mathbf{P}\left( S_{n-j}\leq ta_{k}+\sqrt{ta_{k}}%
,L_{n-j}\geq 0\right) \mathbf{P}_{0}^{+}\left( \mathcal{A}%
_{u.s}|Z_{0}=i\right)
\end{eqnarray*}%
as $n\rightarrow \infty $. Hence, using properties of regular varying
functions and (\ref{BasicAsymptotic}), it is not difficult to deduce from (%
\ref{vaf2}) that
\begin{eqnarray}
&&\lim_{n\gg k\rightarrow \infty }\sup_{-\sqrt{ta_{k}}\leq q\leq w}\frac{%
\mathbf{P}\left( S_{n-j}\leq ta_{k}-q,L_{n-j}\geq 0\right) }{\mathbf{P}%
\left( S_{n-j}\leq ta_{k},L_{n-j}\geq 0\right) }  \notag \\
&=&\lim_{n\gg k\rightarrow \infty }\frac{\mathbf{P}\left( S_{n-j}\leq ta_{k}+%
\sqrt{ta_{k}},L_{n-j}\geq 0\right) }{\mathbf{P}\left( S_{n-j}\leq
ta_{k},L_{n-j}\geq 0\right) }=1.  \label{d2}
\end{eqnarray}%
This estimate shows that, for all $n\geq j$
\begin{eqnarray*}
&&\sup_{-\sqrt{ta_{k}}\leq q\leq w}\frac{\mathbf{E}\left[ H_{n-j}(i);%
\mathcal{T}_{r-j,n-j}(ya_{m});S_{n-j}\leq ta_{k}-q,L_{n-j}\geq 0\right] }{%
\mathbf{P}\left( S_{n-j}\leq ta_{k},L_{n-j}\geq 0\right) } \\
&&\qquad \leq \sup_{-\sqrt{ta_{k}}\leq q\leq w}\frac{\mathbf{P}\left(
S_{n-j}\leq ta_{k}-q,L_{n-j}\geq 0\right) }{\mathbf{P}\left( S_{n-j}\leq
ta_{k},L_{n-j}\geq 0\right) } \\
&&\qquad \quad =\frac{\mathbf{P}\left( S_{n-j}\leq ta_{k}+\sqrt{ta_{k}}%
,L_{n-j}\geq 0\right) }{\mathbf{P}\left( S_{n-j}\leq ta_{k},L_{n-j}\geq
0\right) }\leq C.
\end{eqnarray*}%
Setting $H_{n}=H_{n-j}(i)$ in (\ref{Cond3}) \ and using (\ref{ConvAs})--(\ref%
{d2}), we obtain by Theorem~\ref{T_Minim}, that, for each $q\in \left[ -%
\sqrt{ta_{k}},w\right] $
\begin{eqnarray}
&&\lim_{n\gg k\gg m\rightarrow \infty }\frac{\mathbf{E}\left[ H_{n-j}(i);%
\mathcal{T}_{r-j,n-j}(ya_{m}),S_{n-j}\leq ta_{k}-q,L_{n-j}\geq 0\right] }{%
\mathbf{P}\left( S_{n-j}\leq ta_{k},L_{n-j}\geq 0\right) }  \notag \\
&&\quad =\lim_{n\gg k\gg m\rightarrow \infty }\mathbf{E}\left[ H_{n-j}(i)|%
\mathcal{T}_{r-j,n-j}(ya_{m}),S_{n-j}\leq ta_{k}-q,L_{n-j}\geq 0\right]
\notag \\
&&\qquad \qquad \quad \times \lim_{n\gg k\gg m\rightarrow \infty }\mathbf{P}%
\left( \mathcal{T}_{r-j,n-j}(ya_{m})|S_{n-j}\leq ta_{k}-q,L_{n-j}\geq
0\right)  \notag \\
&&\qquad \qquad \quad \times \lim_{n\gg k\rightarrow \infty }\frac{\mathbf{P}%
\left( S_{n-j}\leq ta_{k}-q,L_{n-j}\geq 0\right) }{\mathbf{P}\left(
S_{n-j}\leq ta_{k},L_{n-j}\geq 0\right) }  \notag \\
&&\qquad \qquad \quad =\mathbf{E}_{0}^{+}\left[ H_{\infty }(i)\right] Q(y)=%
\mathbf{P}_{0}^{+}(\mathcal{A}_{u.s}|Z_{0}=i)Q(y).  \label{Vott1}
\end{eqnarray}%
Further, for $j\in \mathbb{N}_{0}$ and $K\in \mathbb{N}\cup \left\{ \infty
\right\} $ we set
\begin{equation}
\Theta _{w,j}(K)=\Theta _{j}(K):=\sum_{i=1}^{K}\mathbf{P}_{w}\left(
Z_{j}=i,\tau _{j}=j\right) \mathbf{P}_{0}^{+}\left( \mathcal{A}%
_{u.s}|Z_{0}=i\right) .  \label{Vott2}
\end{equation}%
Note that
\begin{equation*}
\Theta _{j}(\infty )=\sum_{i=1}^{\infty }\mathbf{P}_{w}\left( Z_{j}=i,\tau
_{j}=j\right) \mathbf{P}_{0}^{+}\left( \mathcal{A}_{u.s}|Z_{0}=i\right)
<\infty
\end{equation*}%
in view of (\ref{DefTheta}).

\bigskip Applying the dominated convergence theorem we deduce from (\ref%
{Vott1}) and (\ref{Vott2}) that for fixed $K\in \mathbb{N}$
\begin{eqnarray}
&&\lim_{n\gg k\gg m\rightarrow \infty }\frac{\mathbf{P}_{w}\left( \mathcal{T}%
_{r,n}(ya_{m}),S_{n}\leq ta_{k},Z_{n}>0,Z_{j}\leq K,\tau _{n}=j,S_{j}>-\sqrt{%
ta_{k}}\right) }{\mathbf{P}\left( S_{n-j}\leq ta_{k},L_{n-j}\geq 0\right) }
\notag \\
&&\qquad =\lim_{n\gg k\gg m\rightarrow \infty }\int_{-\sqrt{ta_{k}}%
}^{w}\sum_{i=1}^{K}\mathbf{P}_{w}\left( S_{j}\in dq,Z_{j}=i,\tau
_{j}=j\right)  \notag \\
&&\qquad \qquad \times \frac{\mathbf{E}\left[ H_{n-j}(i);\mathcal{T}%
_{r-j,n-j}(ya_{m}),S_{n-j}\leq ta_{k}-q,L_{n-j}\geq 0\right] }{\mathbf{P}%
\left( S_{n-j}\leq ta_{k},L_{n-j}\geq 0\right) }  \notag \\
&&\qquad \qquad =Q(y)\int_{-\infty }^{w}\sum_{i=1}^{K}\mathbf{P}_{w}\left(
S_{j}\in dq,Z_{j}=i,\tau _{j}=j\right) \mathbf{P}_{0}^{+}\left( \mathcal{A}%
_{u.s}|Z_{0}=i\right)  \notag \\
&&\qquad \qquad \qquad =Q(y)\Theta _{j}(K).  \label{Dominate}
\end{eqnarray}%
Therefore,
\begin{eqnarray*}
&&\lim_{n\gg k\gg m\rightarrow \infty }\frac{\mathbf{P}_{w}\left( \mathcal{T}%
_{r,n}(ya_{m}),S_{n}\leq ta_{k},Z_{n}>0,Z_{j}\leq K,\tau _{n}\leq J,S_{j}>-%
\sqrt{ta_{k}}\right) }{\Theta \mathbf{P}\left( S_{n-j}\leq
ta_{k},L_{n-j}\geq 0\right) } \\
&&\qquad \qquad \qquad =\frac{Q(y)}{\Theta }\sum_{j=0}^{J}\Theta _{j}(K).
\end{eqnarray*}%
Hence, letting $K\rightarrow \infty $ and $J\rightarrow \infty $ \ and
taking into account the definition of $\Theta $ in (\ref{DefTheta}) and (\ref%
{AsympMain1}), we get (\ref{Afan0}), which implies
\begin{equation*}
\lim_{n\gg k\gg m\rightarrow \infty }\frac{\mathbf{P}_{w}\left( \mathcal{T}%
_{r,n}(ya_{m}),S_{n}\leq ta_{k},Z_{n}>0\right) }{\mathbf{P}_{w}\left(
S_{n}\leq ta_{k},Z_{n}>0\right) }=Q(y),
\end{equation*}%
that completes the proof of Lemma \ref{L_SmallSurvival}.

\begin{lemma}
\label{L_IntermSurvival} If Conditions B1 and B2 are valid and $\min
(k,r,m)\rightarrow \infty $ as $n\rightarrow \infty ,$ then, for any $w>0,$ $%
y\geq 0,$ and $t>0$
\begin{equation*}
\mathbf{P}_{w}\left( S_{\tau _{r,n}}\leq ya_{k\bot r}|S_{n}\leq
ta_{k},Z_{n}>0\right) \sim \hat{Q}(y),
\end{equation*}%
where $\hat{Q}(y)$ is specified by (\ref{DefQ^}).
\end{lemma}

\textbf{Proof}. Setting $\mathcal{T}_{r,n}(x)=\left\{ S_{\tau _{r,n}}\leq
x\right\} $ in Lemma \ref{L_preliminary}, we see that, similarly to the
proof of Lemma \ref{L_SmallSurvival}, it is necessary to analyse, for fixed $%
j\in \lbrack 0,J],J<n,$ and $K\in \mathbb{N}$ the asymptotic behavior of the
probability
\begin{eqnarray*}
&&\mathbf{P}_{w}\left( \mathcal{T}_{r,n}(ya_{k\bot r}),S_{n}\leq
ta_{k},Z_{n}>0,Z_{j}\leq K,\tau _{n}=j,S_{j}>-\sqrt{ta_{k}}\right) \\
&&\quad =\int_{-\sqrt{ta_{k\bot r}}}^{w}\sum_{i=1}^{K}\mathbf{P}_{w}\left(
S_{j}\in dq,Z_{j}=i,\tau _{j}=j\right) \\
&&\quad \quad \times \mathbf{E}\left[ H_{n-j}(i);\mathcal{T}%
_{r-j,n-j}(ya_{k\bot r}-q),S_{n-j}\leq ta_{k}-q,L_{n-j}\geq 0\right]
\end{eqnarray*}%
if $n\gg k\rightarrow \infty $ and other parameters vary according to the
ranges described in (\ref{DefQ^}).

Using Theorem \ref{T_beginningNew}, relations (\ref{Cond22}), (\ref{DefQ^}),
and (\ref{d2}), we conclude that, for each $q\in \left[ -\sqrt{ta_{k\bot r}}%
,w\right] $ and $y>0$
\begin{eqnarray*}
&&\frac{\mathbf{E}\left[ H_{n-j}(i);\mathcal{T}_{r-j,n-j}(ya_{k\bot
r}-q),S_{n-j}\leq ta_{k}-q,L_{n-j}\geq 0\right] }{\mathbf{P}\left(
S_{n-j}\leq ta_{k},L_{n-j}\geq 0\right) } \\
&=&\mathbf{E}\left[ H_{n-j}(i)|\mathcal{T}_{r-j,n-j}(ya_{k\bot
r}-q),S_{n-j}\leq ta_{k}-q,L_{n-j}\geq 0\right] \\
&&\qquad \times \mathbf{P}\left( \mathcal{T}_{r-j,n-j}(ya_{k\bot
r}-q)|S_{n-j}\leq ta_{k}-q,L_{n-j}\geq 0\right) \\
&&\qquad \times \frac{\mathbf{P}\left( S_{n-j}\leq ta_{k}-q,L_{n-j}\geq
0\right) }{\mathbf{P}\left( S_{n-j}\leq ta_{k},L_{n-j}\geq 0\right) } \\
&&\qquad \qquad \sim \mathbf{E}_{0}^{+}\left[ H_{\infty }(i)\right] \hat{Q}%
(y)=\mathbf{P}_{0}^{+}(\mathcal{A}_{u.s}|Z_{0}=i)\hat{Q}(y).
\end{eqnarray*}

Repeating now with evident changes the part of the proof of Lemma \ref%
{L_SmallSurvival} starting from formula (\ref{Dominate}), using $\mathcal{T}%
_{r-j,n-j}(ya_{k\bot r}-q)$ for $\mathcal{T}_{r-j,n-j}(ya_{m})$ and
replacing $Q(y)$ by $\hat{Q}(y)$ we get
\begin{equation*}
\lim_{n\gg k\rightarrow \infty }\frac{\mathbf{P}_{w}\left( S_{\tau
_{r,n}}\leq ya_{k\bot r},S_{n}\leq ta_{k},Z_{n}>0\right) }{\mathbf{P}%
_{w}\left( S_{n}\leq ta_{k},Z_{n}>0\right) }=\hat{Q}(y),
\end{equation*}%
as required.

Lemma \ref{L_IntermSurvival} is proved.

For integers $0\leq r\leq n$ we introduce the rescaled process $%
X^{r,n}=\left\{ X_{u}^{r,n},0\leq u\leq 1\right\} $, given by
\begin{equation*}
X_{u}^{r,n}:=e^{-S_{r+\left[ (n-r)u\right] }}Z_{r+\left[ (n-r)u\right]
},\quad 0\leq u\leq 1,\quad
\end{equation*}%
and describing the dynamics of the generation sizes of a branching process
in random environment.

\begin{lemma}
\label{L_Skorokh} (Theorem 1 in \cite{VD2023B}). Let Conditions B1 and B2 be
valid and let $r_{1},r_{2},...,$ be a sequence of positive integers such
that $r_{n}\leq n$ and $r_{n}\rightarrow \infty $. If $\varphi (n),n\in
\mathbb{N},$ is a deterministic function such that $a_{n}\gg \varphi
(n)\rightarrow \infty $ as $n\rightarrow \infty $, then
\begin{equation*}
\mathcal{L}\left( X_{u}^{r_{n},n},0\leq u\leq 1|S_{n}\leq \varphi
(n),Z_{n}>0\right) \underset{n\rightarrow \infty }{\Longrightarrow }\mathcal{%
L}\left( W_{u},0\leq u\leq 1\right) ,
\end{equation*}%
where $\left\{ W_{u},0\leq u\leq 1\right\} $ is a stochastic process with
a.s. constant paths, that is
\begin{equation*}
\mathbf{P}\left( W_{u}=W\text{ for all }u\in (0,1]\right) =1
\end{equation*}%
for some random variable $W$. Furthermore,%
\begin{equation*}
\mathbf{P}\left( 0<W<\infty \right) =1.
\end{equation*}
\end{lemma}

\section{Reduced processes in random environment \label{Sec5}}

We now come back to the reduced processes in i.i.d. random environment and
label particles of the $r$th generation by numbers $1,2,...,Z_{r}$ in an
arbitrary but fixed way. Denote by $Z_{r,n}^{(q)},\,q=1,2,\ldots ,Z_{r}$ the
offspring size at moment $n$ of the population generated by the $q$th
particle among those existed at moment $0\leq r<n$. Observe, that%
\begin{eqnarray*}
Z_{r,n} &=&\sum_{q=1}^{Z_{r}}I\left\{ Z_{r,n}^{(q)}>0\right\} \\
&=&Z_{r}\left( 1-F_{r,n}(0)\right) +\sum_{q=1}^{Z_{r}}\left[ I\left\{
Z_{r,n}^{(q)}>0\right\} -\left( 1-F_{r,n}(0)\right) \right] .
\end{eqnarray*}

Our first goal is to show that given $\left\{ S_{n}\leq
ta_{k},Z_{n}>0\right\} $ the random variable $Z_{r,n}$ may be well
approximated by the product $Z_{r}\left( 1-F_{r,n}(0)\right) $ for a wide
range of parameter $r=r(n)$ as $n\rightarrow \infty $. To this aim we
introduce, for $N>0$ and $r<n$ the event

\begin{equation*}
\mathcal{K}_{r,n}(N):=\left\{ \left\vert \sum_{q=1}^{Z_{r}}\left[
I(Z_{r,n}^{(q)}>0)-(1-F_{r,n}(0))\right] \right\vert >\sqrt{%
NZ_{r}(1-F_{r,n}(0))}\right\} ,
\end{equation*}%
denote by $\mathcal{\bar{K}}_{r,n}(N)$ the complementary event to $\mathcal{K%
}_{r,n}(N)$ and prove the following statement.

\begin{lemma}
\label{L_submain} Let Conditions B1 and B2 be valid. Then, for any sequence $%
r=r(n)\rightarrow \infty $ as $n\rightarrow \infty $ such that $r\left(
n\right) <n$
\begin{equation*}
\lim_{N\rightarrow \infty }\limsup_{n\gg k\rightarrow \infty }\mathbf{P}(%
\mathcal{K}_{r,n}(N)|S_{n}\leq ta_{k},Z_{n}>0)=0.
\end{equation*}
\end{lemma}

\textbf{Proof.} Using asymptotic relation (\ref{AsympMain1}) and limit (\ref%
{Prelim1}) with $\mathcal{T}_{n}=\mathcal{K}_{r,n}(N),$ we have
\begin{eqnarray}
&&\lim_{J\rightarrow \infty }\limsup_{n\gg k\rightarrow \infty }\frac{%
\mathbf{P}(\mathcal{K}_{r,n}(N);S_{n}\leq ta_{k},Z_{n}>0,\tau _{n}>J)}{%
\mathbf{P}(S_{n}\leq ta_{k},Z_{n}>0)}  \notag \\
&&\quad =\lim_{J\rightarrow \infty }\limsup_{n\gg k\rightarrow \infty }\frac{%
\mathbf{P}(\mathcal{K}_{r,n}(N);S_{n}\leq ta_{k},Z_{n}>0,\tau _{n}>J)}{%
\Theta \mathbf{P}\left( S_{n}\leq ta_{k},L_{n}\geq 0\right) }=0.
\label{FirstTerm}
\end{eqnarray}%
In this and following lemmas we sometimes write for brevity $\mathbf{P}_{%
\mathcal{F}}(\bullet )$, $\mathbf{E}_{\mathcal{F}}\left[ \bullet \right] ,$
and $\mathbf{D}_{\mathcal{F}}\left[ \bullet \right] $ for the conditional
probabilities, expectations and variances calculated for the $\sigma $%
-algebra $\mathcal{F}_{r,n}$, generated by the random variables $\left\{
Z_{i},i=0,1,...,r\right\} $ and $F_{1}(s),F_{2}(s),\ldots ,F_{n}(s)$.

Clearly, for each $q\in \mathbb{N}$
\begin{equation*}
\mathbf{E}_{\mathcal{F}}\left[ I(Z_{r,n}^{(q)}>0)\right] =1-F_{r,n}(0),\quad
\mathbf{D}_{\mathcal{F}}\left[ I(Z_{r,n}^{(q)}>0)\right]
=(1-F_{r,n}(0))F_{r,n}(0).
\end{equation*}%
Using these relations, observing that, for any $J<r$
\begin{equation*}
\mathbf{P}(Z_{r}>0|\mathcal{E})\leq \min_{0\leq j\leq J}\mathbf{P}(Z_{j}>0|%
\mathcal{E})\leq \min_{0\leq j\leq J}\mathbf{E}\left[ Z_{j}|\mathcal{E}%
\right] =\min_{0\leq j\leq J}e^{S_{j}}=e^{S_{\tau _{J}}},
\end{equation*}%
and applying Chebyshev's inequality to evaluate the probability under the
expectation sign \ in the following formula we obtain
\begin{eqnarray*}
\mathbf{P(}\mathcal{K}_{r,n}(N),S_{n} &\leq &ta_{k},Z_{n}>0,\tau _{n}\leq J)
\\
&\leq &\mathbf{P(}\mathcal{K}_{r,n}(N),S_{n}\leq ta_{k},Z_{r}>0,\tau
_{n}\leq J) \\
&=&\mathbf{E}\left[ I\left\{ \mathcal{K}_{r,n}(N)\right\} I\left\{
Z_{r}>0,S_{n}\leq ta_{k},\tau _{n}\leq J\right\} \right] \\
&=&\mathbf{E}\left[ I\left\{ Z_{r}>0,S_{n}\leq ta_{k},\tau _{n}\leq
J\right\} \mathbf{P}_{\mathcal{F}}\left( \mathcal{K}_{r,n}(N)\right) \right]
\\
&\leq &\mathbf{E}\left[ \frac{I\left\{ Z_{r}>0,S_{n}\leq ta_{k},\tau
_{n}\leq J\right\} }{NZ_{r}(1-F_{r,n}(0))}\mathbf{D}_{\mathcal{F}}\left[
\sum_{q=1}^{Z_{r}}I(Z_{r,n}^{(q)}>0)\right] \right] \\
&\leq &N^{-1}\mathbf{E}\left[ I\{Z_{r}>0\}I\{S_{n}\leq ta_{k},\tau _{n}\leq
J\}\right] \\
&=&N^{-1}\mathbf{E}\left[ \mathbf{P}\left( Z_{r}>0|\mathcal{E}\right)
I\{S_{n}\leq ta_{k},\tau _{n}\leq J\}\right] \\
&\leq &N^{-1}\mathbf{E}\left[ e^{S_{\tau _{J}}};S_{n}\leq ta_{k},\tau
_{n}\leq J\right] \leq N^{-1}\mathbf{E}\left[ e^{S_{\tau _{n}}};S_{n}\leq
ta_{k}\right] \\
&\leq &CN^{-1}\mathbf{P}\left( S_{n}\leq ta_{k},L_{n}\geq 0\right) ,
\end{eqnarray*}%
where, to justify the last inequality, we have used the estimate
\begin{equation*}
\mathbf{E}\left[ e^{S_{\tau _{n}}};S_{n}\leq ta_{k}\right] \sim C_{1}\mathbf{%
P}\left( S_{n}\leq ta_{k},L_{n}\geq 0\right) ,C_{1}\in (0,\infty ),
\end{equation*}%
established in Theorem 1 in \cite{VDD2024} for the case $n\gg k\rightarrow
\infty $.

Thus, for each fixed $J$
\begin{equation*}
\lim_{N\rightarrow \infty }\limsup_{n\gg k\rightarrow \infty }\frac{\mathbf{P%
}\left( \mathcal{K}_{r,n}(N),S_{n}\leq ta_{k},Z_{n}>0,\tau _{n}\leq J\right)
}{\mathbf{P}\left( S_{n}\leq ta_{k},L_{n}\geq 0\right) }=0.
\end{equation*}%
This estimate combined with (\ref{AsympMain1}) and (\ref{FirstTerm}) proves
the lemma. \bigskip

Lemma \ref{L_submain} and the definition of $\mathcal{K}_{r,n}(N)$ give the
following algorithm for fining the conditional distribution of $\log Z_{r,n}$
given $\left\{ S_{n}\leq ta_{k},Z_{n}>0\right\} $ and $n\gg k\rightarrow
\infty :$

First, it is necessary to show that if the event $\left\{ S_{n}\leq
ta_{k},Z_{n}>0\right\} $ has occurred, then the product $Z_{r}(1-F_{r,n}(0))$
\ tends to infinity with probability 1 as $n\gg k\rightarrow \infty $. If
this is the case, then the following representation is valid on the set $%
\mathcal{\bar{K}}_{r,n}(N):$
\begin{equation*}
\log Z_{r,n}=\log \left[ Z_{r}(1-F_{r,n}(0))(1+\varepsilon _{r,n}(k))\right]
,
\end{equation*}%
where $\varepsilon _{r,n}(k)\rightarrow 0$ with probability 1 as $n\gg
k\rightarrow \infty $.

With this result in hands we replace the condition $\left\{ S_{n}\leq
ta_{k},L_{n}\geq 0\right\} $ by the condition $\left\{ S_{n}\leq
ta_{k},Z_{n}>0\right\} $ and, finally, given $\left\{ S_{n}\leq
ta_{k},Z_{n}>0\right\} $ find the distribution of the random variable $\log %
\left[ Z_{r}\left( 1-F_{r,n}(0)\right) \right] $ under an appropriate scaling%
$.$

We proceed to fulfill the program.

For $y>\log 3$ and $r<n$ introduce the events
\begin{equation*}
\mathcal{G}_{r,n}(y)=\{\,\ln (e+Z_{r,n})>y\,\}
\end{equation*}%
and
\begin{equation*}
\mathcal{C}_{r,n}(y)=\left\{ 0<Z_{r}(1-F_{r,n}(0))<1+e^{y}\right\} .
\end{equation*}

\begin{lemma}
\label{L_DeltaNegl0} For any $x>0$ and any sequence $d_{m}=d_{m}(n)%
\rightarrow \infty $ as $n\rightarrow \infty $ we have
\begin{equation}
\lim_{\varepsilon \downarrow 0}\limsup_{n\gg k\rightarrow \infty
}\max_{0\leq r<n}\mathbf{P}(\mathcal{G}_{r,n}(xd_{m})\mathcal{C}%
_{r,n}(\varepsilon d_{m})|S_{n}\leq ta_{k},L_{n}\geq 0)=0.  \label{Est_C}
\end{equation}
\end{lemma}

\textbf{Proof}. Using Chebyshev's inequality we have
\begin{eqnarray}
&&\mathbf{P}(\mathcal{G}_{r,n}(xd_{m})\mathcal{C}_{r,n}(\varepsilon
d_{m}),S_{n}\leq ta_{k},L_{n}\geq 0)  \notag \\
&&\quad =\mathbf{E}\left[ \mathbf{P}_{\mathcal{F}}(\mathcal{G}_{r,n}(xd_{m})%
\mathcal{C}_{r,n}(\varepsilon d_{m}),S_{n}\leq ta_{k},L_{n}\geq 0)\right]
\notag \\
&&\quad =\mathbf{E}\left[ I\left\{ \mathcal{C}_{r,n}(\varepsilon
d_{m}),S_{n}\leq ta_{k},L_{n}\geq 0\right\} \mathbf{P}_{\mathcal{F}}(%
\mathcal{G}_{r,n}(xd_{m}))\right]  \notag \\
&&\quad \leq \frac{1}{x^{2}d_{m}^{2}}\mathbf{E}\Big[\,I\left\{ \mathcal{C}%
_{r,n}(\varepsilon d_{m}),S_{n}\leq ta_{k},L_{n}\geq 0\right\} \mathbf{E}_{%
\mathcal{F}}[\,\ln ^{2}(e+Z_{r,n})\,]\Big].  \label{Est_C1}
\end{eqnarray}

Since
\begin{equation*}
(\ln ^{2}(e+z))^{\prime \prime }=\frac{2}{(e+z)^{2}}(1-\ln (e+z))<0,\quad
z>0,
\end{equation*}%
the function $\ln ^{2}(e+z)$ is concave on the set $\left\{ z>0\right\} $.
This fact allows us to apply Jensen's inequality to the inner expectation on
the right-hand side of (\ref{Est_C1}) and obtain the estimate
\begin{equation*}
\mathbf{E}_{\mathcal{F}}\left[ \ln ^{2}(e+Z_{r,n})\right] \leq \ln ^{2}(e+%
\mathbf{E}_{\mathcal{F}}\left[ Z_{r,n}\right] )=\ln
^{2}(e+Z_{r}(1-F_{r,n}(0))).
\end{equation*}%
Therefore,
\begin{eqnarray*}
&&\mathbf{E}\Big[\,I\left\{ \mathcal{C}_{r,n}(\varepsilon d_{m}),S_{n}\leq
ta_{k},L_{n}\geq 0\right\} \mathbf{E}_{\mathcal{F}}[\,\ln ^{2}(e+Z_{r,n})\,]%
\Big] \\
&\leq &\ln ^{2}(e+e^{\varepsilon d_{m}})\mathbf{P}\left( S_{n}\leq
ta_{k},L_{n}\geq 0\right) \leq \ln ^{2}e^{2\varepsilon d_{m}}\mathbf{P}%
\left( S_{n}\leq ta_{k},L_{n}\geq 0\right) \\
&\leq &4\varepsilon ^{2}d_{m}^{2}\mathbf{P}\left( S_{n}\leq ta_{k},L_{n}\geq
0\right)
\end{eqnarray*}%
if $\varepsilon d_{m}\geq \log 3$, that yields the inequality
\begin{equation*}
\mathbf{P}(\mathcal{G}_{r,n}(xd_{m})\mathcal{C}_{r,n}(\varepsilon
d_{m}),S_{n}\leq ta_{k},L_{n}\geq 0)\leq \frac{4\varepsilon ^{2}}{x^{2}}%
\mathbf{P}\left( S_{n}\leq ta_{k},L_{n}\geq 0\right)
\end{equation*}%
which proves the lemma.

The following statement differs from Lemma \ref{L_DeltaNegl0} by one point
only: we use instead of the condition $\left\{ S_{n}\leq ta_{k},L_{n}\geq
0\right\} $ the condition $\left\{ S_{n}\leq ta_{k},Z_{n}>0\right\} .$

\begin{lemma}
\label{L_REplacement0} Let Conditions $B1$ and $B2$ be valid and $n\gg
k\rightarrow \infty .$ For any $x>0$ and any sequence $d_{m}=d_{m}(n)%
\rightarrow \infty $ as $n\rightarrow \infty $ we have
\begin{equation}
\lim_{\varepsilon \downarrow 0}\limsup_{n\gg k\rightarrow \infty
}\max_{0\leq r<n}\mathbf{P}(\mathcal{G}_{r,n}(xd_{m})\mathcal{C}%
_{r,n}(\varepsilon d_{m})|S_{n}\leq ta_{k},Z_{n}>0)=0.  \label{Est_C_replace}
\end{equation}
\end{lemma}

\textbf{Proof}. It follows from Lemma \ref{L_preliminary} with
\begin{equation*}
\mathcal{T}_{n}=\mathcal{T}_{r,n}=\mathcal{G}_{r,n}(xd_{m})\mathcal{C}%
_{r,n}(\varepsilon d_{m})
\end{equation*}%
that to prove Lemma \ref{L_REplacement0} it is sufficient to show that
\begin{equation}
\lim_{\varepsilon \downarrow 0}\lim_{J\rightarrow \infty }\lim_{K\rightarrow
\infty }\limsup_{n\gg k\rightarrow \infty }\max_{0\leq r<n}\frac{\mathbf{P}%
\left( \mathcal{T}_{r,n},S_{n}\leq ta_{k},Z_{j}\leq K,\tau _{n}\leq
J,S_{\tau _{n}}>-\sqrt{ta_{k}}\right) }{\mathbf{P}\left( S_{n}\leq
ta_{k},Z_{n}>0\right) }=0.  \label{Main_term}
\end{equation}%
To demonstrate the validity of this statement we fix $j\in \lbrack 0,J],J<n,$
and $K\geq 1$ and write
\begin{eqnarray*}
&&\mathbf{P}\left( \mathcal{T}_{r,n},S_{n}\leq ta_{k},Z_{j}\leq K,\tau
_{n}=j,S_{j}>-\sqrt{ta_{k}}\right) \\
&&\qquad \qquad =\int_{-\sqrt{ta_{k}}}^{0}\sum_{i=1}^{K}\mathbf{P}\left(
S_{j}\in dq,Z_{j}=i,\tau _{j}=j\right) \\
&&\qquad \qquad \times \mathbf{E}\left[ \mathbf{P}\left( \mathcal{T}%
_{r-j,n-j}|\mathcal{E},Z_{0}=i\right) ;S_{n-j}\leq ta_{k}-q,L_{n-j}\geq 0%
\right] .
\end{eqnarray*}%
Now using the same estimates as in the proof of Lemma \ref{L_DeltaNegl0}, we
get
\begin{eqnarray*}
&&\mathbf{P}\left( \mathcal{T}_{r-j,n-j}|\mathcal{E},Z_{0}=i\right) \\
&=&\mathbf{P}\left( \ln
(e+Z_{r-j,n-j})>xd_{m},0<Z_{r-j}(1-F_{r-j,n-j}(0))<1+e^{\varepsilon d_{m}}|%
\mathcal{E},Z_{0}=i\right) \\
&\leq &\frac{1}{x^{2}d_{m}^{2}}\mathbf{E}\Big[\,I(C_{r-j,n-j}(\varepsilon
d_{m}))\mathbf{E}_{\mathcal{F}}\left[ \,\ln ^{2}(e+Z_{r-j,n-j})\right] |\,%
\mathcal{E},Z_{0}=i]\Big] \\
&\leq &\frac{1}{x^{2}d_{m}^{2}}\mathbf{E}\Big[\,I(C_{r-j,n-j}(\varepsilon
d_{m}))\,\ln ^{2}(e+\mathbf{E}_{\mathcal{F}}\left[ Z_{r-j,n-j}\right] )|\,%
\mathcal{E},Z_{0}=i\Big] \\
&=&\frac{1}{x^{2}d_{m}^{2}}\mathbf{E}\Big[\,I(C_{r-j,n-j}(\varepsilon
d_{m}))\,\ln ^{2}(e+Z_{r-j}(1-F_{r-j,n-j}(0)))|\,\mathcal{E},Z_{0}=i]\Big] \\
&\leq &\frac{4\varepsilon ^{2}}{x^{2}}.
\end{eqnarray*}%
Thus,
\begin{eqnarray*}
&&\mathbf{P}\left( \mathcal{T}_{r,n},S_{n}\leq ta_{k},Z_{j}\leq K,\tau
_{n}=j,S_{j}>-\sqrt{ta_{k}}\right) \\
&\leq &\frac{4\varepsilon ^{2}}{x^{2}}\int_{-\sqrt{ta_{k}}}^{0}\sum_{i=1}^{K}%
\mathbf{P}\left( S_{j}\in dq,Z_{j}=i,\tau _{j}=j\right) \mathbf{P}\left(
S_{n-j}\leq ta_{k}-q,L_{n-j}\geq 0\right) .
\end{eqnarray*}%
Observing that, for fixed $j\leq J$
\begin{eqnarray*}
&&\lim_{n\gg k\rightarrow \infty }\sup_{q\in (-\sqrt{ta_{k}},0]}\frac{%
\mathbf{P}\left( S_{n-j}\leq ta_{k}-q,L_{n-j}\geq 0\right) }{\mathbf{P}%
\left( S_{n}\leq ta_{k},L_{n}\geq 0\right) } \\
&&\qquad \leq \lim_{n\gg k\rightarrow \infty }\frac{\mathbf{P}\left(
S_{n-j}\leq ta_{k}+\sqrt{ta_{k}},L_{n-j}\geq 0\right) }{\mathbf{P}\left(
S_{n}\leq ta_{k},L_{n}\geq 0\right) }=1,
\end{eqnarray*}%
we obtain for each $J$, sufficiently large $n$ and $\Theta $ specified by (%
\ref{DefTheta}),
\begin{eqnarray*}
&&\frac{\mathbf{P}\left( \mathcal{T}_{r,n},S_{n}\leq ta_{k},Z_{j}\leq K,\tau
_{n}\leq J,S_{\tau _{n}}>-\sqrt{ta_{k}}\right) }{\mathbf{P}\left( S_{n}\leq
ta_{k},Z_{n}>0\right) } \\
&\leq &\frac{4\varepsilon ^{2}}{x^{2}}\sum_{j=0}^{J}\int_{-\sqrt{ta_{k}}}^{0}%
\mathbf{P}\left( S_{j}\in dq,1\leq Z_{j}\leq K,\tau _{j}=j\right) \frac{%
\mathbf{P}\left( S_{n-j}\leq ta_{k}-q,L_{n-j}\geq 0\right) }{\mathbf{P}%
\left( S_{n}\leq ta_{k},Z_{n}>0\right) } \\
&\leq &\frac{8\varepsilon ^{2}}{x^{2}\Theta }\sum_{j=0}^{J}\int_{-\sqrt{%
ta_{k}}}^{0}\mathbf{P}\left( S_{j}\in dq,1\leq Z_{j}\leq K,\tau
_{j}=j\right) \frac{\mathbf{P}\left( S_{n-j}\leq ta_{k}+\sqrt{ta_{k}}%
,L_{n-j}\geq 0\right) }{\mathbf{P}\left( S_{n}\leq ta_{k},L_{n}\geq 0\right)
} \\
&\leq &\frac{16\varepsilon ^{2}}{x^{2}\Theta }\sum_{j=0}^{J}\int_{-\sqrt{%
ta_{k}}}^{0}\mathbf{P}\left( S_{j}\in dq,1\leq Z_{j}\leq K,\tau _{j}=j\right)
\\
&\leq &\frac{16\varepsilon ^{2}}{x^{2}\Theta }\sum_{j=0}^{\infty }\mathbf{P}%
\left( Z_{j}>0,\tau _{j}=j\right) \leq \frac{C\varepsilon ^{2}}{x^{2}\Theta }%
,
\end{eqnarray*}%
where for the last step we have used (\ref{Sparr}). This evidently yields (%
\ref{Main_term}).

Lemma \ref{L_REplacement0} is proved.

\subsection{Proof of Theorem \textbf{\protect\ref{T_reducedSmall}}}

We find in this section the limit of the conditional probability
\begin{equation*}
\mathbf{P}\left( \log Z_{r,n}-S_{r}\leq za_{m}|\mathcal{R}\left(
ta_{k},n\right) \right)
\end{equation*}%
as $n\gg k\gg m=n-r\rightarrow \infty $.

Our first step in the analysis of the probability we are interested in is
the following important claim concerning the distribution of the minimum of
the associated random walk $\left\{ S_{i},i=0,1,2,...\right\} $ on the
interval $[r,n]$ given $n\gg k\gg m=n-r\rightarrow \infty $.

\begin{lemma}
\label{L_Difference_small} If Conditions B1 and B2 are valid then, for any $%
y\geq 0$
\begin{equation*}
\lim_{n\gg k\gg m\rightarrow \infty }\mathbf{P}\left( \frac{1}{a_{k}}S_{\tau
_{r,n}}>y\Big|\mathcal{R}\left( ta_{k},n\right) \right) =1-\left( \frac{%
t\wedge y}{t}\right) ^{\alpha \rho +1}.
\end{equation*}
\end{lemma}

\textbf{Proof}. We write the decomposition
\begin{equation*}
S_{\tau _{r,n}}=S_{\tau _{r,n}}-S_{r}+S_{r}-S_{n}+S_{n}.
\end{equation*}%
Note that the assumption $k\gg m$ implies $a_{k}\gg a_{m}$. Keeping this
relation in mind we conclude by Lemma \ref{L_SmallSurvival} that, for any $%
\varepsilon >0$
\begin{equation}
\lim_{n\gg k\gg m\rightarrow \infty }\mathbf{P}\left( \frac{\left\vert
S_{\tau _{r,n}}-S_{r}\right\vert }{a_{k}}>\varepsilon \Big|\mathcal{R}\left(
ta_{k},n\right) \right) =0.  \label{Remm1}
\end{equation}%
Further, it follows from Lemma \ref{L_Skorokh} that given Conditions B1 and
B2
\begin{equation}
\lim_{n\gg k\gg m\rightarrow \infty }\mathcal{L}\left( \frac{Z_{r}}{e^{S_{r}}%
}\Big|\mathcal{R}\left( ta_{k},n\right) \right) =\mathcal{L}\left( W\right) ,
\label{Mart_ratio}
\end{equation}%
where the random variable $W$ possesses the property
\begin{equation*}
\mathbf{P}\left( 0<W<\infty \right) =1.
\end{equation*}%
Besides, according to point 3) of Theorem 4 in \cite{VDD2023}, for any $z\in
(-\infty ,+\infty )$
\begin{eqnarray*}
&&\lim_{n\gg k\gg m\rightarrow \infty }\mathbf{P}\left( \frac{1}{a_{m}}(\log
Z_{r}-S_{n})\leq z\Big|\mathcal{R}\left( ta_{k},n\right) \right) \\
&=&\lim_{n\gg k\gg m\rightarrow \infty }\mathbf{P}\left( \frac{1}{a_{m}}%
\left( \log \frac{Z_{r}}{e^{S_{r}}}+S_{r}-S_{n}\right) \leq z\Big|\mathcal{R}%
\left( ta_{k},n\right) \right) =\mathbf{P}(Y_{1}\leq z).
\end{eqnarray*}%
This relation combined with (\ref{Mart_ratio}) implies
\begin{equation*}
\lim_{n\gg k\gg m\rightarrow \infty }\mathbf{P}\left( \frac{1}{a_{m}}%
(S_{r}-S_{n})\leq z\Big|\mathcal{R}\left( ta_{k},n\right) \right) =\mathbf{P}%
\left( Y_{1}\leq z\right) .
\end{equation*}%
As a result, for any $\varepsilon >0$
\begin{equation}
\lim_{n\gg k\gg m\rightarrow \infty }\mathbf{P}\left( \frac{\left\vert
S_{r}-S_{n}\right\vert }{a_{k}}>\varepsilon \Big|\mathcal{R}\left(
ta_{k},n\right) \right) =0.  \label{Remm2}
\end{equation}%
Finally, for any $y\geq 0$
\begin{equation}
\lim_{n\gg k\rightarrow \infty }\mathbf{P}\left( \frac{1}{a_{k}}S_{n}\geq y%
\Big|\mathcal{R}\left( ta_{k},n\right) \right) =1-\left( \frac{t\wedge y}{t}%
\right) ^{\alpha \rho +1}  \label{Remm3}
\end{equation}%
in view of (\ref{AsympMain1}), (\ref{Regular1}), and (\ref{AsympV}).
Combining (\ref{Remm1})--(\ref{Remm3}) proves the lemma.

Introduce the notation
\begin{equation*}
\eta _{q}:=\frac{F_{q}^{\prime \prime }(1)}{\left( F_{q}^{\prime }(1)\right)
^{2}},\quad q=1,2,....,
\end{equation*}%
and set
\begin{equation*}
O_{r,n}:=1+\sum_{q=r}^{n-1}\eta _{q+1}e^{S_{\tau _{r,n}}-S_{q}}\text{ and }%
O_{n}:=O_{0,n}.
\end{equation*}

\begin{lemma}
\label{L_OConveregnce} If Conditions B1 and B2 are valid, then
\begin{equation*}
\lim_{N\rightarrow \infty }\limsup_{n\rightarrow \infty }\mathbf{P}\left(
O_{n}>2N\right) =0.
\end{equation*}
\end{lemma}

\textbf{Proof}. Clearly,
\begin{eqnarray*}
\mathbf{P}\left( O_{n}>2N\right) &=&\sum_{j=0}^{n}\mathbf{P}\left(
1+\sum_{q=0}^{n-1}\eta _{q+1}e^{S_{j}-S_{q}}>2N,\tau _{n}=j\right) \\
&\leq &\Xi _{1}(n,N)+\Xi _{2}(n,N),
\end{eqnarray*}%
where
\begin{eqnarray*}
\Xi _{1}(n,N):= &&\sum_{j=0}^{n}\mathbf{P}\left( 1+\sum_{q=0}^{j-1}\eta
_{q+1}e^{S_{j}-S_{q}}>N,\tau _{n}=j\right) ; \\
\Xi _{2}(n,N):= &&\sum_{j=0}^{n}\mathbf{P}\left( \sum_{q=j}^{n-1}\eta
_{q+1}e^{S_{j}-S_{q}}>N,\tau _{n}=j\right) .
\end{eqnarray*}%
Set
\begin{equation*}
S_{l}^{\prime }:=S_{j+l}-S_{j},\quad \eta _{l}^{\prime }:=\eta _{j+l},\quad
l=0,1,...,n-j,\quad L_{n-j}^{\prime }:=\min_{0\leq l\leq n-j}S_{l}^{\prime }.
\end{equation*}%
We have
\begin{align*}
\mathbf{P}\left( \sum_{q=j}^{n-1}\eta _{q+1}e^{S_{j}-S_{q}}>N,\tau
_{n}=j\right) & =\mathbf{P}\left( \tau _{j}=j\right) \mathbf{P}\left(
\sum_{l=0}^{n-j-1}\eta _{l+1}^{\prime }e^{-S_{l}^{\prime
}}>N,L_{n-j}^{\prime }\geq 0\right) \\
& =\mathbf{P}\left( M_{j}<0\right) \mathbf{P}\left( L_{n-j}\geq 0\right)
\mathbf{P}\left( H_{n-j-1}>N|L_{n-j}\geq 0\right) \\
& =\mathbf{P}\left( \tau _{n}=j\right) \mathbf{P}\left(
H_{n-j-1}>N|L_{n-j}\geq 0\right) ,
\end{align*}%
where%
\begin{equation*}
H_{n-j-1}:=\sum_{l=0}^{n-j-1}\eta _{l+1}e^{-S_{l}}.
\end{equation*}%
Clearly,
\begin{equation*}
\lim_{n\rightarrow \infty }H_{n-j-1}=H_{\infty }:=\sum_{l=0}^{\infty }\eta
_{l+1}e^{-S_{l}}\quad \text{a.s.},
\end{equation*}%
where
\begin{equation*}
H_{\infty }<\infty \quad \mathbf{P}_{0}^{+}\text{- a.s.}
\end{equation*}%
according to Lemma 2.7 in \cite{agkv}. This estimate and Lemma 2.5 in \cite%
{agkv}, applied to random variable $e^{-\lambda H_{n-j-1}}$ with $\lambda
>0, $ imply
\begin{equation*}
\lim_{N\rightarrow \infty }\lim_{n\rightarrow \infty }\mathbf{P}\left(
H_{n-j-1}>N|L_{n-j}\geq 0\right) =\lim_{N\rightarrow \infty }\mathbf{P}%
_{0}^{+}\left( H_{\infty }>N\right) =0.
\end{equation*}%
Since, for each fixed $J\in \mathbb{N}$
\begin{eqnarray*}
&&\lim_{N\rightarrow \infty }\max_{0\leq n-j\leq J}\mathbf{P}\left(
\sum_{q=j}^{n-1}\eta _{q+1}e^{S_{j}-S_{q}}>N,\tau _{n}=j\right) \\
&&\qquad \leq \lim_{N\rightarrow \infty }\max_{0\leq n-j\leq J}\mathbf{P}%
\left( \sum_{q=j}^{n-1}\eta _{q+1}>N,\tau _{n}=j\right) \\
&&\qquad \leq \lim_{N\rightarrow \infty }\max_{0\leq n-j\leq J}\mathbf{P}%
\left( \sum_{q=0}^{J}\eta _{q+1}>N\right) =0,
\end{eqnarray*}%
we have
\begin{eqnarray}
&&\lim_{N\rightarrow \infty }\lim_{n\rightarrow \infty }\Xi _{2}(n,N)  \notag
\\
&&=\lim_{N\rightarrow \infty }\lim_{n\rightarrow \infty }\sum_{j=0}^{n}%
\mathbf{P}\left( \tau _{n}=j\right) \mathbf{P}\left( H_{n-j-1}>N|L_{n-j}\geq
0\right) =0.  \label{Thet1}
\end{eqnarray}

We now set
\begin{equation*}
S_{j-q}^{\prime }:=S_{j}-S_{q},\quad \eta _{j-q}^{\prime }=\eta _{q+1},\quad
q=0,1,...,j;\quad M_{j}^{\prime }:=\max_{1\leq i\leq j}S_{i}^{\prime }.
\end{equation*}%
Then
\begin{eqnarray*}
&&\mathbf{P}\left( 1+\sum_{q=0}^{j-1}\eta _{q+1}e^{S_{j}-S_{q}}>N,\tau
_{n}=j\right) \\
&&\qquad =\mathbf{P}\left( 1+\sum_{q=1}^{j}\eta _{q}^{\prime
}e^{S_{q}^{\prime }}>N,M_{j}^{\prime }<0\right) \mathbf{P}\left( L_{n-j}\geq
0\right) \\
&&\qquad =\mathbf{P}\left( M_{j}<0\right) \mathbf{P}\left( L_{n-j}\geq
0\right) \mathbf{P}\left( U_{j}>N|M_{j}<0\right) ,
\end{eqnarray*}%
where
\begin{equation*}
U_{j}:=1+\sum_{q=1}^{j}\eta _{q}e^{S_{q}}\rightarrow U_{\infty
}:=1+\sum_{q=1}^{\infty }\eta _{q}e^{S_{q}}\quad \text{a.s.}
\end{equation*}%
as $j\rightarrow \infty $. According to Lemma 3.1 in \cite{ABGV2011}
\begin{equation*}
U_{\infty }<\infty \quad \mathbf{P}_{0}^{-}\text{- a.s.,}
\end{equation*}%
where the measure $\mathbf{P}_{0}^{-}$ is specified in (\ref{StayNegative}).
This fact and Lemma 2.3 in \cite{ABGV2013}, applied to $e^{-\lambda U_{j}}$
with $\lambda >0,$ imply
\begin{equation*}
\lim_{N\rightarrow \infty }\lim_{n\rightarrow \infty }\mathbf{P}\left(
U_{j}>N|M_{j}<0\right) =\lim_{N\rightarrow \infty }\mathbf{P}_{0}^{-}\left(
U_{\infty }>N\right) =0.
\end{equation*}

Now, by the same arguments as we have used to establish (\ref{Thet1}), it is
not difficult to conclude that
\begin{equation*}
\lim_{N\rightarrow \infty }\lim_{n\rightarrow \infty }\Xi _{1}(n,N)=0.
\end{equation*}

This relation and (\ref{Thet1}) imply that
\begin{equation*}
\lim_{N\rightarrow \infty }\limsup_{n\rightarrow \infty }\mathbf{P}\left(
O_{n}>2N\right) \leq \lim_{N\rightarrow \infty }\lim_{n\rightarrow \infty
}\left( \Xi _{1}(n,N)+\Xi _{2}(n,N)\right) =0.
\end{equation*}

Lemma \ref{L_OConveregnce} is proved.

Set
\begin{equation*}
\Delta _{r,n}=\log \frac{Z_{r}(1-F_{r,n}(0))}{e^{S_{\tau _{r,n}}}}=\log
\frac{Z_{r}}{e^{S_{r}}}+\log \frac{1-F_{r,n}(0)}{e^{S_{\tau _{r,n}}-S_{r}}}.
\end{equation*}%
The next lemma combining with Lemma \ref{L_Difference_small} shows, in
particular, that if the event $\mathcal{R}\left( ta_{k},n\right) $ has
occurred$,$ then $Z_{r}(1-F_{r,n}(0))\rightarrow \infty $ with probability 1
as $n\gg k\gg n-r\rightarrow \infty $.

\begin{lemma}
\label{L_DeltaNegl} If Conditions B1 and B2 are valid, then, for any fixed $%
\varepsilon >0$
\begin{equation}
\lim_{n\gg k\gg m\rightarrow \infty }\mathbf{P}\left( \frac{1}{a_{m}}%
\left\vert \Delta _{r,n}\right\vert >\varepsilon \Big|\mathcal{R}\left(
ta_{k},n\right) \right) =0.  \label{DeltaNegl}
\end{equation}
\end{lemma}

\textbf{Proof}. It is known (see, for instance, relations (2.2) and (2.3) in
\cite{GK2000}) that, for any $j=1,2,...,$
\begin{equation*}
1-F_{0,j}(0)\geq \left( e^{-S_{j}}+\sum_{q=0}^{j-1}\eta
_{q+1}e^{-S_{q}}\right) ^{-1}=e^{S_{\tau _{j}}}\left( 1+\sum_{q=0}^{j-1}\eta
_{q+1}e^{S_{\tau _{j}}-S_{q}}\right) ^{-1}.
\end{equation*}%
This inequality implies the estimate
\begin{equation*}
\log (1-F_{r,n}(0))\geq \left( S_{\tau _{r,n}}-S_{r}\right) -\log O_{r,n},
\end{equation*}%
which, in view of evident inequalities
\begin{equation*}
1-F_{r,n}(0)\leq 1-F_{r,\tau _{r,n}}(0)\leq e^{S_{\tau _{r,n}}-S_{r}},
\end{equation*}%
shows that
\begin{equation*}
0\geq \log \frac{1-F_{r,n}(0)}{e^{S_{\tau _{r,n}}-S_{r}}}\geq -\log O_{r,n}.
\end{equation*}%
It is easy to check that, for any $\varepsilon >0$
\begin{eqnarray*}
&&\mathbf{P}\left( \log O_{r,n}>\varepsilon a_{m},\mathcal{R}\left(
ta_{k},n\right) \right) \\
&&\quad \leq \mathbf{P}\left( S_{r}-S_{\tau _{r.n}}>\sqrt{a_{m}a_{k}},%
\mathcal{R}\left( ta_{k},n\right) \right) \\
&&\quad +\mathbf{P}\left( \log O_{r,n}>\varepsilon a_{m},S_{r}\leq ta_{k}+%
\sqrt{a_{m}a_{k}},Z_{r}>0\right) .
\end{eqnarray*}%
By Lemma \ref{L_SmallSurvival}
\begin{equation*}
\lim_{n\gg k\gg m\rightarrow \infty }\frac{\mathbf{P}\left( S_{r}-S_{\tau
_{r.n}}>\sqrt{a_{m}a_{k}},\mathcal{R}\left( ta_{k},n\right) \right) }{%
\mathbf{P}\left( \mathcal{R}\left( ta_{k},n\right) \right) }=0.
\end{equation*}%
Now we observe that the random variable $O_{r,n}$ is independent of the
event $\left\{ S_{r}\leq ta_{k}+\sqrt{a_{m}a_{k}},Z_{r}>0\right\} $. Hence,
taking into account (\ref{AsympMain1}), we get
\begin{eqnarray*}
&&\mathbf{P}\left( \log O_{r,n}>\varepsilon a_{m},S_{r}\leq ta_{k}+\sqrt{%
a_{m}a_{k}},Z_{r}>0\right) \\
&&\qquad =\mathbf{P}\left( O_{m}>e^{\varepsilon a_{m}}\right) \mathbf{P}%
\left( S_{r}\leq ta_{k}+\sqrt{a_{m}a_{k}},Z_{r}>0\right) \\
&&\qquad \leq \mathbf{P}\left( O_{m}>e^{\varepsilon a_{m}}\right) \mathbf{P}%
\left( S_{r}\leq (t+1)a_{k},Z_{r}>0\right)
\end{eqnarray*}%
for sufficiently large $m=n-r$ and $n\gg k\gg m$. This estimate and Lemma %
\ref{L_OConveregnce} show that, for any $\varepsilon >0$
\begin{equation}
\lim_{n\gg k\gg m\rightarrow \infty }\mathbf{P}\left( \frac{1}{a_{m}}%
\left\vert \log \frac{1-F_{r,n}(0)}{e^{S_{\tau _{r,n}}-S_{r}}}\right\vert >%
\frac{\varepsilon }{2}\Big|\mathcal{R}\left( ta_{k},n\right) \right) =0.
\label{Delta1}
\end{equation}

Further, in view of (\ref{Mart_ratio}), for any $\varepsilon >0$ and any
sequence $d_{k}\rightarrow \infty $ as $k\rightarrow \infty $
\begin{equation}
\lim_{n\gg k\gg m\rightarrow \infty }\mathbf{P}\left( \frac{1}{d_{k}}%
\left\vert \log \frac{Z_{r}}{e^{S_{r}}}\right\vert >\frac{\varepsilon }{2}%
\Big|\mathcal{R}\left( ta_{k},n\right) \right) =0.  \label{Delta2}
\end{equation}%
Combining (\ref{Delta1}) and (\ref{Delta2}) gives (\ref{DeltaNegl}).

Lemma \ref{L_DeltaNegl} is proved.

\textbf{Proof of Theorem \ref{T_reducedSmall}}. It follows from Lemma \ref%
{L_DeltaNegl}, the continuity of $Q(x)$, $x\in (-\infty ,+\infty ),$ and
Lemma \ref{L_SmallSurvival} that%
\begin{eqnarray}
&&\lim_{n\gg k\gg m\rightarrow \infty }\mathbf{P}\left( \frac{1}{a_{m}}\log %
\left[ \frac{Z_{r}(1-F_{,n}(0))}{e^{S_{r}}}\right] \leq x\Big|\mathcal{R}%
\left( ta_{k},n\right) \right)  \notag \\
&&  \notag \\
&=&\lim_{n\gg k\gg m\rightarrow \infty }\mathbf{P}\left( \frac{\Delta
_{r,n}+S_{\tau _{r,n}}-S_{r}}{a_{m}}\leq x\Big|\mathcal{R}\left(
ta_{k},n\right) \right) =Q(x).  \label{ProductLimit}
\end{eqnarray}%
Further, by Lemma \ref{L_Difference_small}, for any $y>0$
\begin{eqnarray*}
&&\liminf_{n\gg k\gg m\rightarrow \infty }\mathbf{P}\left( \frac{1}{a_{m}}%
\log \left[ Z_{r}(1-F_{r,n}(0))\right] \geq y\Big|\mathcal{R}\left(
ta_{k},n\right) \right) \\
&&\qquad =\liminf_{n\gg k\gg m\rightarrow \infty }\mathbf{P}\left( \frac{1}{%
a_{m}}\Delta _{r,n}+\frac{1}{a_{m}}S_{\tau _{r,n}}\geq y\Big|\mathcal{R}%
\left( ta_{k},n\right) \right) \\
&&\qquad =\liminf_{n\gg k\gg m\rightarrow \infty }\mathbf{P}\left( \frac{1}{%
a_{m}}S_{\tau _{r,n}}\geq y\Big|\mathcal{R}\left( ta_{k},n\right) \right) =1.
\end{eqnarray*}%
In particular, for any $M>0$
\begin{equation*}
\lim_{n\gg k\gg m\rightarrow \infty }\mathbf{P}\left( Z_{r}(1-F_{r,n}(0))>M|%
\mathcal{R}\left( ta_{k},n\right) \right) =1.
\end{equation*}

Now to complete the proof of Theorem \ref{T_reducedSmall} \ it is sufficient
to observe that in view of Lemmas~\ref{L_submain} and \ref{L_DeltaNegl} and
relation (\ref{ProductLimit})
\begin{eqnarray*}
&&\lim_{n\gg k\gg m\rightarrow \infty }\mathbf{P}\left( \frac{1}{a_{m}}\log
\frac{Z_{r,n}}{e^{S_{r}}}\leq x\Big|\mathcal{R}\left( ta_{k},n\right) \right)
\\
&=&\lim_{N\rightarrow \infty }\lim_{n\gg k\gg m\rightarrow \infty }\mathbf{P}%
\left( \frac{1}{a_{m}}\log \frac{Z_{r,n}}{e^{S_{r}}}\leq x,\mathcal{\bar{K}}%
_{l,n}(N)\Big|\mathcal{R}\left( ta_{k},n\right) \right) \\
&=&\lim_{N\rightarrow \infty }\lim_{n\gg k\gg m\rightarrow \infty }\mathbf{P}%
\left( \frac{1}{a_{m}}\log \left[ \frac{Z_{r}(1-F_{r,n}(0))}{e^{S_{r}}}%
\right] \leq x,\mathcal{\bar{K}}_{l,n}(N)\Big|\mathcal{R}\left(
ta_{k},n\right) \right) \\
&=&\lim_{n\gg k\gg m\rightarrow \infty }\mathbf{P}\left( \frac{1}{a_{m}}\log %
\left[ \frac{Z_{r}(1-F_{r,n}(0))}{e^{S_{r}}}\right] \leq x\Big|\mathcal{R}%
\left( ta_{k},n\right) \right) =Q(x).
\end{eqnarray*}

Theorem \ref{T_reducedSmall} is proved\textbf{.}

\subsection{Limit theorem for $n\gg k\sim \protect\theta m$}

In this section we consider the case $n\gg k\sim \theta m=\theta (n-r)$ for
some $\theta >0$ and start by proving the following statement.

\begin{lemma}
\label{L_Difference_Interem} If Condition B1 is valid, then, for any $z\in
(-t,\infty )$%
\begin{equation*}
\lim_{n\gg k\sim \theta m\rightarrow \infty }\mathbf{P}\left(
S_{r}-S_{n}\leq za_{k}|\mathcal{B}(ta_{k},n)\right) =A_{2}(t,\theta ,z),
\end{equation*}%
where%
\begin{equation*}
A_{2}(t,\theta ,z)=\frac{\alpha \rho +1}{t^{\alpha \rho +1}}%
\int_{0}^{t+z}v^{\alpha \rho }\mathbf{P}\left( -v\theta ^{1/\alpha }\leq
\min_{0\leq s\leq 1}Y_{s},-z\theta ^{1/\alpha }\leq Y_{1}\leq \left(
t-v\right) \theta ^{1/\alpha }\right) dv.
\end{equation*}
\end{lemma}

\begin{remark}
\label{Rem2} It is not difficult to see that
\begin{eqnarray*}
A_{2}(t,\theta ,\infty ):= &&\lim_{z\rightarrow \infty }A_{2}(t,\theta ,z) \\
&=&\frac{\alpha \rho +1}{t^{\alpha \rho +1}}\int_{0}^{\infty }v^{\alpha \rho
}\mathbf{P}\left( -v\theta ^{1/\alpha }\leq \min_{0\leq s\leq
1}Y_{s},Y_{1}\leq t\theta ^{1/\alpha }-v\theta ^{1/\alpha }\right) dv \\
&=&\frac{\alpha \rho +1}{\left( t\theta ^{1/\alpha }\right) ^{\alpha \rho +1}%
}\int_{0}^{\infty }w^{\alpha \rho }\mathbf{P}\left( -w\leq \min_{0\leq s\leq
1}Y_{s},Y_{1}\leq t\theta ^{1/\alpha }-w\right) dw \\
&=&A\left( t\theta ^{1/\alpha },t\theta ^{1/\alpha }\right) .
\end{eqnarray*}%
This representation and Remark \ref{Rem1} show that $A_{2}(t,\theta ,\infty
)=1$, and, therefore,
\begin{equation}
\lim_{z\rightarrow \infty }\lim_{n\gg k\sim \theta m\rightarrow \infty }%
\mathbf{P}\left( S_{r}\geq za_{k}|\mathcal{B}(ta_{k},n)\right) =0.
\label{NeglInterm}
\end{equation}
\end{remark}

\textbf{\ Proof of Lemma \ref{L_Difference_Interem}}. Setting
\begin{equation*}
S_{j}^{\prime }:=S_{r+j}-S_{r},\quad j=0,1,...,m,\quad L_{m}^{\prime
}:=\min_{0\leq j\leq m}S_{j}^{\prime },
\end{equation*}%
we have
\begin{eqnarray*}
&&\mathbf{P}\left( S_{r}-S_{n}\leq za_{k},S_{n}\leq ta_{k},L_{n}\geq 0\right)
\\
&=&\int_{0}^{t+z}\mathbf{P}\left( S_{r}\in a_{k}dv,L_{r}\geq 0\right)
\mathbf{P}\left( -a_{k}v\leq L_{m}^{\prime },-S_{m}^{\prime }\leq
za_{k},S_{m}^{\prime }\leq \left( t-v\right) a_{k}\right) \\
&=&\int_{0}^{t+z}\mathbf{P}\left( S_{r}\in a_{k}dv,L_{r}\geq 0\right)
\mathbf{P}\left( -a_{k}v\leq L_{m},-za_{k}\leq S_{m}\leq \left( t-v\right)
a_{k}\right) .
\end{eqnarray*}%
It is known (see formula (5.2) in Theorem 5.1 in \cite{CC2013} ), that
\begin{equation*}
\frac{\mathbf{P}\left( S_{r}\in a_{k}dv,L_{r}\geq 0\right) }{a_{k}dv}\sim
\frac{g_{\alpha ,\beta }(0)}{ra_{r}}V^{+}(a_{k}v)\sim \frac{g_{\alpha ,\beta
}(0)}{na_{n}}V^{+}(a_{k}v)
\end{equation*}%
as $n\sim r\rightarrow \infty $ uniformly in $va_{k}\ll a_{r}$.

\bigskip Moreover, in view of (\ref{Regular1})
\begin{equation*}
\frac{V^{+}(a_{k}v)}{V^{+}(a_{k}v)}\rightarrow v^{\alpha \rho }
\end{equation*}%
as $k\rightarrow \infty $ uniformly in $v\in \lbrack 0,t=z]$. Therefore, as $%
n>>k\rightarrow \infty $%
\begin{eqnarray*}
&&\mathbf{P}\left( S_{r}-S_{n}\leq za_{k},\mathcal{B}(ta_{k},n)\right) \\
&\sim &\frac{g_{\alpha ,\beta }(0)a_{k}}{na_{n}}\int_{0}^{t+z}V^{+}(a_{k}v)%
\mathbf{P}\left( -va_{k}\leq L_{m},-za_{k}\leq S_{m}\leq \left( t-v\right)
a_{k}\right) dv \\
&\sim &\frac{g_{\alpha ,\beta }(0)a_{k}V^{+}(a_{k})}{na_{n}}%
\int_{0}^{t+z}v^{\alpha \rho }\mathbf{P}\left( -v\frac{a_{k}}{a_{m}}\leq
\frac{L_{m}}{a_{m}},-z\frac{a_{k}}{a_{m}}\leq \frac{S_{m}}{a_{m}}\leq \left(
t-v\right) \frac{a_{k}}{a_{m}}\right) dv.
\end{eqnarray*}%
\bigskip

\bigskip On the other hand, if $a_{k}\sim \theta ^{1/\alpha }a_{m},$ then,
for any $v\in \lbrack 0,t+z]$
\begin{eqnarray*}
&&\lim_{k\sim \theta m\rightarrow \infty }\mathbf{P}\left( -v\frac{a_{k}}{%
a_{m}}\leq \frac{L_{m}}{a_{m}},-z\frac{a_{k}}{a_{m}}\leq \frac{S_{m}}{a_{m}}%
\leq \left( t-v\right) \frac{a_{k}}{a_{m}}\right) \\
&&\qquad =\mathbf{P}\left( -v\theta ^{1/\alpha }\leq \min_{0\leq s\leq
1}Y_{s},-z\theta ^{1/\alpha }\leq Y_{1}\leq \left( t-v\right) \theta
^{1/\alpha }\right) .
\end{eqnarray*}%
This relation and the dominated convergence theorem imply that
\begin{eqnarray*}
&&\mathbf{P}\left( S_{r}-S_{n}\leq za_{k},\mathcal{B}(ta_{k},n)\right) \\
&\sim &\frac{g_{\alpha ,\beta }(0)a_{k}V^{+}(a_{k})}{na_{n}}%
\int_{0}^{t+z}v^{\alpha \rho }\mathbf{P}\left( -v\frac{a_{k}}{a_{m}}\leq
\frac{L_{m}}{a_{m}},-z\frac{a_{k}}{a_{m}}\leq \frac{S_{m}}{a_{m}}\leq \left(
t-v\right) \frac{a_{k}}{a_{m}}\right) dv \\
&\sim &\frac{g_{\alpha ,\beta }(0)a_{k}V^{+}(a_{k})}{na_{n}}%
\int_{0}^{t+z}v^{\alpha \rho }\mathbf{P}\left( -v\theta ^{1/\alpha }\leq
\min_{0\leq s\leq 1}Y_{s},-z\theta ^{1/\alpha }\leq Y_{1}\leq \left(
t-v\right) \theta ^{1/\alpha }\right) dv
\end{eqnarray*}%
as $n\gg k\sim \theta m\rightarrow \infty $. Since
\begin{equation*}
\mathbf{P}\left( \mathcal{B}(ta_{k},n)\right) \sim g_{\alpha ,\beta
}(0)b_{n}\int_{0}^{ta_{k}}V^{+}(u)du\sim g_{\alpha ,\beta }(0)b_{n}\frac{%
t^{\alpha \rho +1}a_{k}V^{+}(a_{k})}{\alpha \rho +1}
\end{equation*}%
as $n\gg k\rightarrow \infty $, it follows from the estimates above that
\begin{equation*}
\lim_{n\gg k\sim \theta m\rightarrow \infty }\mathbf{P}\left(
S_{r}-S_{n}\leq za_{k}|\mathcal{B}(ta_{k},n)\right) =A_{2}(t,\theta ,z).
\end{equation*}

\bigskip Lemma \ref{L_Difference_Interem} is proved\textbf{.}

The aim of the next lemma is to justify the replacement of the condition $%
\mathcal{B}(ta_{k},n)$ by the condition $\mathcal{R}\left( ta_{k},n\right) $.

\begin{lemma}
\label{L_Zero_negl} If Conditions B1 and B2 are\ valid, $k\sim \theta
m=\theta (n-r)$ for some $\theta >0,$ then
\begin{equation*}
\lim_{z\rightarrow \infty }\lim_{n\gg k\sim \theta m\rightarrow \infty }%
\mathbf{P}\left( S_{r}>za_{k}|\mathcal{R}(ta_{k},n)\right) =0.
\end{equation*}
\end{lemma}

\textbf{Proof}. Setting $\mathcal{T}_{n}=\left\{ S_{r}>za_{k}\right\} $ in
Lemma \ref{L_preliminary}, we see that, similarly to the proof of Lemma \ref%
{L_SmallSurvival}, it is necessary to analyse, for fixed $j\in \lbrack 0,J]$
and $K\in \mathbb{N}$ the asymptotic behavior of the probability
\begin{eqnarray*}
&&\mathbf{P}\left( \mathcal{T}_{n},S_{n}\leq ta_{k},Z_{n}>0,Z_{j}\leq K,\tau
_{n}=j,S_{j}>-\sqrt{ta_{k}}\right) \\
&&\qquad =\int_{-\sqrt{ta_{k}}}^{0}\sum_{i=1}^{K}\mathbf{P}\left( S_{j}\in
dq,Z_{j}=i,\tau _{j}=j\right) \\
&&\qquad \qquad \times \mathbf{E}\left[ H_{n-j}(i);S_{r-j}>za_{k}-q,S_{n-j}%
\leq ta_{k}-q,L_{n-j}\geq 0\right]
\end{eqnarray*}%
as $n\gg k\rightarrow \infty $. Clearly,
\begin{eqnarray*}
&&\mathbf{E}\left[ H_{n-j}(i);S_{r-j}>za_{k}-q,S_{n-j}\leq
ta_{k}-q,L_{n-j}\geq 0\right] \\
&&\qquad \leq \mathbf{P}\left( S_{r-j}>za_{k}-q,S_{n-j}\leq
ta_{k}-q,L_{n-j}\geq 0\right) \\
&&\qquad \leq \mathbf{P}\left( S_{r-j}>za_{k},S_{n-j}\leq
2ta_{k},L_{n-j}\geq 0\right)
\end{eqnarray*}%
for all sufficiently large $k$ and $q\in \lbrack -\sqrt{ta_{k}},0].$
Therefore, for each fixed $j\in \lbrack 0,J],$ $K\in \mathbb{N}$ and
sufficiently large $k$
\begin{eqnarray*}
&&\mathbf{P}\left( \mathcal{T}_{n},S_{n}\leq ta_{k},Z_{n}>0,Z_{j}\leq K,\tau
_{n}=j,S_{j}>-\sqrt{ta_{k}}\right) \\
&\leq &\mathbf{P}\left( S_{r-j}>za_{k},S_{n-j}\leq 2ta_{k},L_{n-j}\geq
0\right) \\
&&\times \int_{-\sqrt{ta_{k}}}^{0}\sum_{i=1}^{K}\mathbf{P}\left( S_{j}\in
dq,Z_{j}=i,\tau _{j}=j\right) \\
&\leq &\mathbf{P}\left( S_{r-j}>za_{k},S_{n-j}\leq 2ta_{k},L_{n-j}\geq
0\right) \mathbf{P}\left( 0<Z_{j}\leq K,\tau _{j}=j\right) \\
&\leq &\mathbf{P}\left( S_{r-j}>za_{k},S_{n-j}\leq 2ta_{k},L_{n-j}\geq
0\right) \mathbf{E}\left[ e^{S_{j}},\tau _{j}=j\right] .
\end{eqnarray*}%
This estimate implies the statement of the lemma, since
\begin{equation*}
\sum_{j=0}^{\infty }\mathbf{E}\left[ e^{S_{j}},\tau _{j}=j\right] <\infty
\end{equation*}%
in view of (\ref{ExpFunct}), and, at the same time,
\begin{eqnarray*}
&&\lim_{z\rightarrow \infty }\limsup_{n\gg k\sim \theta m\rightarrow \infty }%
\frac{\mathbf{P}\left( S_{r-j}>za_{k},S_{n-j}\leq 2ta_{k},L_{n-j}\geq
0\right) }{\mathbf{P}\left( S_{n}\leq ta_{k},L_{n}\geq 0\right) } \\
&&\quad =\lim_{z\rightarrow \infty }\limsup_{n\gg k\sim \theta m\rightarrow
\infty }\frac{\mathbf{P}\left( S_{r-j}>za_{k},S_{n-j}\leq
2ta_{k},L_{n-j}\geq 0\right) }{\mathbf{P}\left( S_{n-j}\leq
2ta_{k},L_{n-j}\geq 0\right) } \\
&&\qquad \quad \times \limsup_{n\gg k\rightarrow \infty }\frac{\mathbf{P}%
\left( S_{n-j}\leq 2ta_{k},L_{n-j}\geq 0\right) }{\mathbf{P}\left( S_{n}\leq
ta_{k},L_{n}\geq 0\right) }=0
\end{eqnarray*}%
in view of (\ref{NeglInterm}), (\ref{BasicAsymptotic}) and (\ref{Defb}).

\begin{lemma}
\label{L_DeltaNegll2} Let Conditions B1 and B2 be valid. Then, for any fixed
$\varepsilon >0$%
\begin{equation}
\lim_{n\gg k\sim \theta m\rightarrow \infty }\mathbf{P}\left( \frac{1}{a_{m}}%
\left\vert \Delta _{r,n}\right\vert >\varepsilon \Big|\mathcal{R}\left(
ta_{k},n\right) \right) =0.  \label{DeltaNegl22}
\end{equation}
\end{lemma}

\textbf{Proof}. We use arguments similar to those used to prove Lemma \ref%
{L_DeltaNegl}. For any positive $N$ we have
\begin{eqnarray*}
&&\mathbf{P}\left( O_{r,n}>e^{\varepsilon a_{m}},\mathcal{R}\left(
ta_{k},n\right) \right) \leq \mathbf{P}\left( S_{r}>Na_{k},\mathcal{R}\left(
ta_{k},n\right) \right) \\
&&\qquad \qquad \qquad +\mathbf{P}\left( O_{r,n}>e^{\varepsilon
a_{m}},S_{r}\leq \left( t+N\right) a_{k},Z_{r}>0\right) .
\end{eqnarray*}%
It follows from Lemma \ref{L_Zero_negl} that
\begin{equation*}
\lim_{N\rightarrow \infty }\lim_{n\gg k\sim \theta m\rightarrow \infty }%
\frac{\mathbf{P}\left( S_{r}>Na_{k},\mathcal{R}\left( ta_{k},n\right)
\right) }{\mathbf{P}\left( \mathcal{R}\left( ta_{k},n\right) \right) }=0.
\end{equation*}%
Taking into account (\ref{AsympMain1}), the independence of $O_{r,n}$ of the
event $\left\{ S_{r}\leq ta_{k}+Na_{m},Z_{r}>0\right\} ,$ and Lemma \ref%
{L_OConveregnce} we conclude that, for any $\delta >0$
\begin{eqnarray*}
&&\mathbf{P}\left( O_{r,n}>e^{\varepsilon a_{m}},S_{r}\leq
ta_{k}+Na_{m},Z_{r}>0\right) \\
&&\qquad =\mathbf{P}\left( O_{m}>e^{\varepsilon a_{m}}\right) \mathbf{P}%
\left( S_{r}\leq ta_{k}+Na_{m},Z_{r}>0\right) \\
&&\qquad \quad \leq \delta \mathbf{P}\left( \mathcal{R}\left(
ta_{k},n\right) \right)
\end{eqnarray*}%
for all sufficiently large $n$ and $m$. Therefore, for any $\varepsilon >0$
\begin{equation*}
\lim_{n\gg k\sim \theta m\rightarrow \infty }\mathbf{P}\left( \frac{1}{a_{m}}%
\left\vert \log \frac{1-F_{r,n}(0)}{e^{S_{\tau _{r,n}}-S_{r}}}\right\vert >%
\frac{\varepsilon }{2}\Big|\mathcal{R}\left( ta_{k},n\right) \right) =0.
\end{equation*}

Using now Lemma \ref{L_Skorokh}, we see that, for any $\varepsilon >0$
\begin{equation*}
\lim_{n\gg k\sim \theta m\rightarrow \infty }\mathbf{P}\left( \frac{1}{a_{m}}%
\left\vert \log \frac{Z_{r}}{e^{S_{r}}}\right\vert >\frac{\varepsilon }{2}%
\Big|\mathcal{R}\left( ta_{k},n\right) \right) =0.
\end{equation*}

Combining the obtained estimates proves (\ref{DeltaNegl22}).

The following lemma plays a key role in proving Theorem \ref{T_reducedInterm}%
.

\begin{lemma}
\label{L_intermediate} Let Conditions B1 and B2 be valid. If $k\sim $ $%
\theta m=\theta (n-r)$ for some $\theta >0$ then$,$ for $y\geq 0$
\begin{equation*}
\lim_{n\gg k\rightarrow \infty }\mathbf{P}\left( \frac{1}{a_{m}}\log \left[
Z_{r}(1-F_{r,n}(0))\right] \leq y\Big|\mathcal{R}\left( ta_{k},n\right)
\right) =A(\theta ^{1/\alpha }t,\theta ^{1/\alpha }t\wedge y).
\end{equation*}
\end{lemma}

\textbf{Proof}. Using Lemma \ref{L_DeltaNegll2} and recalling the statement
of Lemma \ref{L_IntermSurvival} we have
\begin{eqnarray*}
&&\mathbf{P}\left( \frac{1}{a_{m}}\log \left[ Z_{r}(1-F_{r,n}(0))\right]
\leq y\Big|\mathcal{R}\left( ta_{k},n\right) \right) \\
&=&\mathbf{P}\left( \frac{1}{a_{m}}\Delta _{r,n}+\frac{1}{a_{m}}S_{\tau
_{r,n}}\leq y\Big|\mathcal{R}\left( ta_{k},n\right) \right) \sim \mathbf{P}%
\left( \frac{1}{a_{m}}S_{\tau _{r,n}}\leq y\Big|\mathcal{R}\left(
ta_{k},n\right) \right) \\
&=&\mathbf{P}\left( S_{\tau _{r,n}}\leq ya_{k}\frac{a_{m}}{a_{k}}\Big|%
\mathcal{R}\left( ta_{k},n\right) \right) \sim A(t\theta ^{1/\alpha },\theta
^{1/\alpha }t\wedge y)
\end{eqnarray*}%
as $n\gg k\sim \theta m\rightarrow \infty $, as desired.

\textbf{Proof of Theorem \ref{T_reducedInterm}}. In view of Lemma \ref%
{L_submain}\ to prove the theorem it is sufficient to show that, for any $%
y\in \lbrack 0,\theta ^{1/\alpha }t]$
\begin{equation}
\lim_{N\rightarrow \infty }\lim_{n\rightarrow \infty }\mathbf{P}\left( \log
Z_{r,n}\leq ya_{m};\mathcal{\bar{K}}_{l,n}(N)|\mathcal{R}\left(
ta_{k},n\right) \right) =A(\theta ^{1/\alpha }t,y).  \label{Check1}
\end{equation}%
Since $Z_{r}\left( 1-F_{r,n}(0)\right) \rightarrow \infty $ with probability
1 as $n\gg m=n-r\rightarrow \infty $ by Lemma \ref{L_intermediate}, we
conclude that (\ref{Check1}) will be proved if we show that
\begin{equation*}
\lim_{N\rightarrow \infty }\lim_{n\rightarrow \infty }\mathbf{P}\left( \log %
\left[ Z_{r}\left( 1-F_{r,n}(0)\right) \right] \leq ya_{m};\mathcal{\bar{K}}%
_{l,n}(N)|\mathcal{R}\left( ta_{k},n\right) \right) =A(\theta ^{1/\alpha
}t,y).
\end{equation*}%
However this fact is an obvious conclusion of the just mentioned Lemmas~\ref%
{L_submain} and \ref{L_intermediate}.

Theorem \ref{T_reducedInterm} is proved.

\subsection{Proof of Theorem \protect\ref{T_GeneralFrac}}

Similarly to the proofs of Theorems \ref{T_reducedSmall} and \ref%
{T_reducedInterm}, the first step in the proof of this theorem is to
demonstrate that given conditions of the theorem
\begin{equation}
\lim_{k\wedge r\rightarrow \infty }\mathbf{P}\left( \frac{1}{a_{k\wedge r}}%
\left\vert \Delta _{r,n}\right\vert >\varepsilon \Big|\mathcal{R}%
_{n}(ta_{k})\right) =0  \label{FirstStep}
\end{equation}%
for any $\varepsilon >0.$ It follows from Lemma \ref{L_Skorokh} that
\begin{equation*}
\lim_{k\wedge r\rightarrow \infty }\mathbf{P}\left( \frac{1}{a_{k\wedge r}}%
\left\vert \log \frac{Z_{r}}{e^{S_{r}}}\right\vert >\frac{\varepsilon }{2}%
\Big|\mathcal{R}_{n}(ta_{k})\right) =0.
\end{equation*}%
Thus, to prove (\ref{FirstStep}) we need to check that
\begin{equation*}
\lim_{k\wedge r\rightarrow \infty }\mathbf{P}\left( \frac{1}{a_{k\wedge r}}%
\left\vert \log \frac{1-F_{r,n}(0)}{e^{S_{\tau _{r,n}}-S_{r}}}\right\vert >%
\frac{\varepsilon }{2}\Big|\mathcal{R}_{n}(ta_{k})\right) =0\;
\end{equation*}%
under the conditions of Theorem \ref{T_GeneralFrac}. This is, indeed, the
case, since it follows from the conditions of Theorem \ref{T_GeneralFrac}
that
\begin{equation*}
\frac{1}{a_{k\wedge r}}\left\vert \log \frac{1-F_{r,n}(0)}{e^{S_{\tau
_{r,n}}-S_{r}}}\right\vert \leq \frac{1}{a_{k\wedge r}}\log O_{r,n}\leq
\frac{1}{a_{k\wedge r}}\log C(n-r)\rightarrow 0\text{ a.s.}
\end{equation*}%
as $k\wedge r\rightarrow \infty $. Since $a_{k\wedge r}\leq a_{k\bot r},$ we
deduce that (\ref{FirstStep}) and Lemma \ref{L_IntermSurvival} imply the
estimate
\begin{eqnarray*}
&&\lim \mathbf{P}\left( \frac{1}{a_{k\bot r}}\log \left[ Z_{r}(1-F_{r,n}(0))%
\right] \leq y\Big|\mathcal{R}\left( ta_{k},n\right) \right) \\
&&\quad =\lim \mathbf{P}\left( \frac{1}{a_{k\bot r}}\left( \Delta
_{r,n}+S_{\tau _{r,n}}\right) \leq y\Big|\mathcal{R}\left( ta_{k},n\right)
\right) \\
&&\qquad =\lim \mathbf{P}\left( \frac{1}{a_{k\bot r}}S_{\tau _{r,n}}\leq y%
\Big|\mathcal{R}\left( ta_{k},n\right) \right) =\hat{Q}(y),
\end{eqnarray*}%
where the limits are calculated when the of parameters $r,k$ and $n$ vary in
accordance with points 1)--3) of Theorem \ref{T_GeneralFrac}. Therefore, if
the event $\mathcal{R}\left( ta_{k},n\right) $\ has occurred and the
parameters $r,k$ and $n$ vary in accordance to the assumptions of Theorem %
\ref{T_GeneralFrac}, then
\begin{equation*}
Z_{r}(1-F_{r,n}(0))\rightarrow \infty \text{ a.s.}
\end{equation*}%
This fact and Lemma \ref{L_submain} \ allow to claim that if Condition (\ref%
{Cond_log}) is valid, then, for any $y\geq 0$

\begin{eqnarray*}
&&\lim_{n\gg k\gg r\rightarrow \infty }\mathbf{P}\left( \log Z_{r,n}\leq
ya_{r}|\mathcal{R}\left( ta_{k},n\right) \right) \\
&&\qquad \qquad =\lim_{N\rightarrow \infty }\lim_{n\gg k\gg r\rightarrow
\infty }\mathbf{P}\left( \log Z_{r,n}\leq ya_{r};\mathcal{\bar{K}}_{r,n}(N)|%
\mathcal{R}\left( ta_{k},n\right) \right) \\
&&\qquad \qquad =\lim_{n\gg k\gg r\rightarrow \infty }\mathbf{P}\left( \frac{%
1}{a_{r}}\log \left[ Z_{r}(1-F_{r,n}(0))\right] \leq y\Big|\mathcal{R}\left(
ta_{k},n\right) \right) \\
&&\qquad \qquad \qquad =C^{\ast \ast }\mathcal{H}(y),
\end{eqnarray*}%
and by the same reasonings, if $n\gg k\sim \theta r$ for some $\theta >0,$
then, for any $y\geq 0$
\begin{equation*}
\lim_{k\rightarrow \infty }\mathbf{P}\left( \log Z_{r,n}\leq ya_{r}|\mathcal{%
R}\left( ta_{k},n\right) \right) =W(\theta ^{1/\alpha }t,\theta ^{1/\alpha
}t\wedge y),
\end{equation*}%
and if $\min \left( r,n-r\right) \gg k$, then, for any $y\geq 0$
\begin{equation*}
\lim_{k\rightarrow \infty }\mathbf{P}\left( \log Z_{r,n}\leq ya_{k}|\mathcal{%
R}\left( ta_{k},n\right) \right) =1-\left( 1-\frac{t\wedge y}{t}\right)
^{\alpha \rho +1}.
\end{equation*}

Theorem \ref{T_GeneralFrac} is proved.

\end{document}